\documentclass[12pt]{amsart}
\usepackage{amsfonts,amsmath,amssymb,amscd,amsthm,enumerate}
\usepackage{amsrefs}
\usepackage{geometry}
\usepackage[all]{xy}
\usepackage{stmaryrd}
\newtheorem{theorem}{Theorem}
\newtheorem{prop}[theorem]{Proposition}
\newtheorem{lem}[theorem]{Lemma}
\newtheorem{cor}[theorem]{Corollary}
\theoremstyle{definition}

\newtheorem{rem}[theorem]{Remark}
\newtheorem{mydef}[theorem]{Definition}
\newtheorem{example}[theorem]{Example}

\renewcommand{\epsilon}{\varepsilon}

\def\ad{{\rm ad}}
\def\Ad{{\rm Ad}}
\def\C{\mathbb{C}}
\def\Z{\mathbb{Z}}
\def\L{\mathbb{L}}
\def\R{\mathbb{R}}

\def\T{\mathbb{T}}

\def\<{\langle}
\def\>{\rangle}

\def\EE{{\mathcal E}}
\def\DD{{\mathcal D}}

\def\ad{{\rm ad}}
\def\Hom{{\rm Hom}}
\def\pr{{\rm pr}}
\begin{document}
\title{Generalized Almost Product Structures and Generalized CRF-structures}
\author{Marco Aldi and Daniele Grandini}
\begin{abstract}
We give several equivalent characterizations of orthogonal subbundles of the generalized tangent bundle defined, up to B-field transform, by almost product and local product structures. We also introduce a pure spinor formalism for generalized CRF-structure and investigate the resulting decomposition of the de Rham operator. As applications we give a characterization of generalized complex manifolds that are locally the product of generalized complex factors and discuss infinitesimal deformations of generalized CRF-structures.
\end{abstract}

\maketitle

\section{Introduction}

Generalized CRF-structures were introduced in \cite{vaismanCRF} as Courant involutive, (not necessarily maximal) isotropic subbundles $L$ of the complexified generalized tangent bundle with no non-trivial totally real section. In this paper we continue the work initiated in \cite{AG2} and focus on generalized CRF-structures 
$L$ such that $L\oplus \overline L=E\otimes \C$, where $E$ is a {\it split structure} i.e.\ a subbundle of the generalized tangent bundle with the property that the restriction of the tautological inner product to $E$ is non-degenerate of signature $(k,k)$. If this is the case we say that $L$ is a generalized CRF-structure on the split structure $E$. Since generalized complex structures \cite{G} and strongly integrable generalized contact structures \cite{PW} are all examples of generalized CRF-structures on split structures, their study is important in order to develop a unified understanding of geometric structures on the generalized tangent bundle.

The goal of this paper is to investigate the geometry of generalized CRF-structures on a particular class of split structures called {\it generalized almost product structures}. By definition a generalized almost product structure is a split structure $E$ whose projection $\pi(E)$ onto the tangent bundle is ``minimal'' in the sense that its rank is half the rank of $E$. As we show, this notion is equivalent to requiring that $(\pi(E),\pi(E^\perp))$ is a classical almost product structure (as defined in \cite{gugenheim-spencer}) or, alternatively, to the condition that the cotangent bundle is a direct sum of its intersections with $E$ and $E^\perp$. This last characterization implies that each generalized almost product structure gives rise to a canonical bigrading on the exterior algebra of differential forms which is a refinement of its standard $\Z$-grading. In particular we obtain a decomposition of the de Rham operator $d=d_E+d_{E^\perp}$, with $d_E$ of bidegree $(1,0)+(-1,2)$. Interestingly, restricting the standard Dorfman bracket to $E$ and composing with the orthogonal projection of the generalized tangent bundle onto $E$ gives rise to a binary operation $\llbracket\, ,\,\rrbracket_E$ which coincides with the derived bracket for $d_E$. We show that $\llbracket\,,\,\rrbracket_E$ (together with the restrictions of the tautological inner product and anchor map to $E$) defines a structure of Courant algebroid on an almost product $E$ if and only if $\pi(E)$ is a foliation. Moreover, $\pi(E)$ and $\pi(E^\perp)$ are complementary foliations if and only if $d_E^2=0$.

An appealing feature of generalized CRF-structures on generalized almost product structures is that they admit an alternate description in terms of pure spinors, which generalizes the spinorial approach to generalized complex structures and  generalized contact structures discussed in \cite{G}, \cite{hitchin03}, \cite{cavalcanti} and\cite{AG}. In particular we show that each generalized CRF-structure gives rise to a canonical (up to shift) $\Z$-grading on complex differential forms. With respect to this grading, $d_{E^\perp}$ is of degree $0$ while $d_E$ decomposes into a components of degree $1$ and $-1$ which in the case of generalized complex structures respectively to the $\partial$ and and $\overline\partial$ operators defined in \cite{G}. 

We also discuss a weaker integrability condition in which $d_E$ is still required to decompose into components of degree $\pm 1$, but no assumption is made on the degree of $d_{E^\perp}$. These more general structures, which we refer to as {\it weak generalized CRF-structures}, contain interesting examples (e.g.\ classical contact structures) that dare not generalized CRF-structures.

In the particular case of generalized CRF-structures on an almost product structure $E$ such that $\pi(E^\perp)$ is a foliation, $d_E$ restricts to a differential on basic forms for this foliation. The grading induced by the generalized CRF-structure and the resulting decomposition of $d_E$ can be restricted to basic forms. In particular, the spinorial approach to transverse generalized complex structures of \cite{wade10} and the basic $dd^{\mathcal J}$-lemma discussed in \cite{razny16} fit naturally into the framework of the present paper. Moreover if $\pi(E)$ is also a foliation, we show that the results of \cite{cavalcanti} on the $\partial\overline\partial$-lemma and the canonical spectral sequence apply to all forms, not just those that are basic with respect to $\pi(E^\perp)$.

In addition to illustrating with examples that our framework effectively unifies previous spinorial approaches to generalized geometry, we offer two applications. The first is a characterization of generalized complex manifolds $(M,J)$ that are locally the product of two generalized complex manifolds in terms of certain integrability conditions satisfied by the restriction of $J$ to an almost product structure. Our second application is a characterization of infinitesimal deformations of (weak) generalized CRF-structures along the lines of the Kodaira-Spencer formalism developed for generalized complex structures in \cite{G}, \cite{li05} and \cite{tomasiello}. We prove that that deformations of weak generalized CRF-structures are governed by a single equation which specializes to the well-known Kodaira-Spencer/Maurer-Cartan equation in the case of generalized complex structures. On the other hand a second equation, stating that the operator that represents the deformation commutes with $d_{E^\perp}$ is required to characterized deformations of generalized CRF-structures.

The paper is organized as follows. In Section 2 we collect some basic facts about $\R$-linear operators acting on differential forms, the preferred language of this paper. In particular, we view sections of the generalized tangent bundle as operators acting on forms in such a way that (up to scaling by a factor of $2$), the tautological inner product coincides with the obvious graded commutators of operators. We also introduce generalized Lie derivatives as well as derived brackets for operators that are not necessarily the de Rham operator, as this level of generality is useful in the bulk of the paper. In Section 3 we introduce generalized almost product structures. After proving several equivalent characterizations of generalized almost product structures among all split structures, we describe the decomposition of the de Rham operator that they induce and the corresponding derived brackets. We then proceed to investigate the additional structure that emerges if one additionally assumes that $\pi(E)$ and/or $\pi(E^\perp)$ is a foliation. In this context, we also introduce our slight generalization (accounting for a possible B-field transform) of the standard notions of basic differential forms and basic complex attached to a foliation. In Section 4 is devoted to Vaisman's generalized $F$-structures, which we view as operators acting on forms. In the particular case in which the kernel of the generalized $F$-structure is an almost product structure (or, more generally, is equipped with a decomposition into isotropic subbundles), we construct a canonical $\Z$-grading on complex differential forms. After illustrating these notions with several examples, we show that with respect to this grading $d_{E^\perp}$ decomposes into components of degree $0$ or $\pm 2$. In Section 5 we investigate the integrability conditions which define (weak) generalized CRF-structures among all generalized $F$-structures. Our characterizations of integrability are intended to be reminiscent of those established for generalized complex structures in \cite{G}, \cite{cavalcanti} and \cite{tomasiello}. In the last part of this section we specialize to the case in which both $\pi(E)$ and $\pi(E^\perp)$ are foliations. In particular, we discuss the role of the $\partial\overline\partial$-lemma in this framework and prove the characterization of local products of generalized complex manifolds mentioned above. The paper ends with Section 6, which is devoted to the study of infinitesimal deformations of (weak) generalized CRF-structures. While (in the spirit of \cite{li05} and \cite{tomasiello}), our results are stated in the language of operators acting on forms, we also remark that in the case in which the image of the generalized $F$-structure is a foliation our finding are in agreement with the standard theory of deformations of Lie bialgebroids developed in \cite{LWX}.

\bigskip\noindent \textbf{Acknowledgments:} Parts of this paper were written while visiting Swarthmore College, the Simons Center for Geometry and Physics and IMPA. We would like to thank these institutions for hospitality and excellent working conditions. We also would like to thank Reimundo Heluani, Ralph Gomez, Janet Talvacchia and Alessandro Tomasiello for inspiring conversations.

\section{Operators on forms}

\begin{mydef}\label{def:1} Unless otherwise specified, we let $M$ be a connected, finite dimensional smooth manifold. We denote by $\Omega_M=\Gamma(\wedge^\bullet T^*M)$ be the graded commutative algebra of $\R$-valued differential forms on $M$. We denote by $\Omega_M^k$, $k=0,\ldots,\dim M$, the graded component with respect to the standard $\Z$-grading and by $\Omega_M^{\bar k}$, $\bar k=\bar 0, \bar 1$ the components of the standard $\Z/2$-grading by parity. We denote by $\EE_M$ the graded algebra of $\R$-linear endomorphisms of $\Omega_M$ and by $\DD_M$ the graded Lie algebra of graded derivations of $\EE_M$. We define the {\it adjoint map} $\ad\in \Hom_\R( \EE_M,\DD_M)$ such that
\begin{equation}\label{eq:1}
\ad_\varphi(\psi)=[\varphi,\psi]=\varphi\circ\psi-(-1)^{kl}\psi\circ\varphi
\end{equation}
for all $\varphi\in \EE_M^k$ and $\psi\in \EE_M^l$.
\end{mydef}

\begin{rem}\label{rem:2}
Left multiplication defines a canonical embedding of $\Omega_M$ into $\EE_M$. On the other hand, if $\varphi\in \EE_M$ is such that $\Omega_M^1\subseteq \ker (\ad_\varphi)$, then $\varphi=\varphi(1)\in \Omega_M$. Therefore, $\varphi\in \Omega_M$ if and only if $\Omega_M^1\subseteq \ker (\ad_\varphi)$. Similarly, $\varphi\in \EE_M$ is $\Omega_M^0$-linear if and only if $\Omega_M^0\subseteq \ker (\ad_\varphi)$. 
\end{rem}

\begin{rem}\label{rem:3}
The canonical embedding of $\Omega_M^1$ in $\EE_M^1$ can be extended to a canonical embedding of sections of the {\it  generalized tangent bundle} $\T M=TM\oplus T^*M$ of $M$ into $\EE_M^1\oplus \EE_M^{-1}$ by identifying each $X\in \Gamma(TM)$ with the interior product $\iota_X$, i.e.\ the unique $\Omega_M^0$-linear operator such that  $[X,\alpha]=\alpha(X)$ for all $\alpha\in \Omega_M^1$. By Remark \ref{rem:2}, $\Gamma(\T M)\subseteq \ker(\ad_\varphi)$ if and only if $\varphi$ is a form on $M$ such that $\Gamma(TM)\subseteq \ker (\ad_\varphi)$, i.e.\ if and only if $\varphi\in \Omega_M^0$.
\end{rem}

\begin{rem}\label{rem:4}
Let $d\in \EE_M^1$ be the de Rham operator on $M$.
Since $T^*M$ is maximal isotropic in $\T M$ (with respect to the $\Omega_M^0$-valued pairing $\langle\,,\,\rangle$ defined as twice the graded commutator) and $\Omega_M^1$ is generated over $\Omega_M^0$ by $d(\Omega_M^0)$, we conclude that
\begin{equation}\label{eq:2}
\Omega_M^1=\{x\in \Gamma(\T M)\,|\,d(\Omega_M^0)\subseteq \ker (\ad_x)\}\,.
\end{equation}
\end{rem}

\begin{mydef}\label{def:5}
We define the {\it generalized Lie derivative associated with $\delta\in \EE_M$} to be $\L^{\!\!^\delta}=\ad\circ \ad_\delta\in \Hom_\R(\EE_M,\DD_M)$ which assigns to each $\varphi\in \EE_M$ the operator $\L_\varphi^{\!\!^\delta}=\ad_{[\varphi,\delta]}$. The {\it derived bracket associated with $\delta$} is $\llbracket\,,\,\rrbracket_\delta:\EE_M\otimes_\R\EE_M\to \EE_M$ defined by $\llbracket \varphi_1,\varphi_2\rrbracket_\delta = \L^{\!\!^\delta}_{\varphi_1} (\varphi_2)$. In the important case of the de Rham operator, we use the shorthand notation $\L=\L^{\!\!^d}$ and $\llbracket \, , \, \rrbracket= \llbracket \, , \, \rrbracket_d$. 
\end{mydef}

\begin{rem}\label{rem:6}
Let $\deg\in\EE_M^0$ be the diagonal operator on $\Omega_M$ with $k$-eigenspace equal to $\Omega^k_M$ for each $k$. Then $\ad_{\deg}$ is diagonal on $\EE_M$ with $l$-eigenspace equal to $\EE_M^l$ for each $l$. Therefore, each derivation  $\delta\in \EE_M$ of degree $l\neq 0$ can be recovered from the corresponding generalized Lie derivative as $\delta=\frac{1}{l}\L^{\!\!^\delta}_{\deg}$. In particular, $d=\L_{\deg}$. 
\end{rem} 

\begin{rem}\label{rem:7}
Using the Jacobi identity for $(\EE_M, [\,,\,])$, conveniently written as
\begin{equation}\label{eq:3}
[\ad_\varphi,\ad_\psi]=\ad_{[\varphi,\psi]}
\end{equation}
for all all $\varphi,\psi\in \EE_M$, it is straightforward to check that the identities
\begin{align}
\ad_{{\llbracket \varphi,\psi \rrbracket}_\delta} &= [\L^{\!\!^\delta}_\varphi, \ad_\psi] \label{eq:4}\\ 
(-1)^{kj}[[\varphi,\psi],\delta] & =  \llbracket \varphi,\psi\rrbracket_\delta  - (-1)^{k(i+j)+ij} \llbracket \psi,\varphi\rrbracket_\delta \label{eq:5}
\end{align}
hold for all $\varphi\in \EE_M^i$, $\psi\in \EE_M^j$ and $\delta\in \EE_M^k$.
\end{rem}

\begin{rem}\label{rem:8}
If $\alpha\in \Omega_M^k$, then $\L_\alpha(f)=(-1)^{k}[f,d\alpha]=0=\ad_\alpha(f)$ for all $f\in \Omega_M^0$. Conversely, if 
\begin{equation}\label{eq:6}
\Omega_M^0 \subseteq \ker (\ad_\varphi)\cap \ker (\L_\varphi)
\end{equation}
then by Remark \ref{rem:2}, $\varphi$ is $\Omega_M^0$-linear and, by \eqref{eq:3}, 
\begin{equation}\label{eq:7}
\ad_\varphi(df)=\ad_\varphi (\ad_d(f))=0
\end{equation}
for all $f\in \Omega_M^0$. Using Remark \ref{rem:2} again we conclude that $\varphi \in \Omega_M$ if and only if \eqref{eq:7} holds. 
\end{rem}

\begin{rem}\label{rem:9}
If $X\in \Gamma(TM)$, then $\L_X$ acts on $\Omega_M$ as the usual Lie derivative. In particular, $\Omega_M^0$ is closed under the action of $\L_\varphi$ for each $\varphi \in \Gamma(TM)\oplus \Omega_M$. Conversely, suppose that $\varphi\in \EE_M$ is $\Omega_M^0$-linear and $\L_\varphi(\Omega_M^0)\subseteq \Omega_M^0$. Then $\L_\varphi$ is a derivation of $\Omega_M^0$ and therefore there exists $X\in \Gamma(TM)$ such that $\Omega_M^0\subseteq \ker(\L_{\varphi-X})$. By Remark \ref{rem:8} we conclude that $\varphi-X\in \Omega_M$ and thus $\varphi\in \Gamma(TM)\oplus \Omega_M$. 
\end{rem}

\begin{mydef}\label{def:10}
If $\varphi\in \EE_M$ is nilpotent, {\it the adjoint automorphism associated with $\varphi$} is $\Ad_\varphi\in {\rm Aut}_{\R}(\EE_M)$ such that $\Ad_\varphi(\psi)=e^\varphi\circ\psi\circ e^{-\varphi}$ for all $\psi\in \EE_M$. 
\end{mydef}

\begin{rem}\label{rem:11}
Since $M$ is finite dimensional, every $\varphi\in \EE_M$ of positive or negative degree is automatically nilpotent. 
\end{rem}

\begin{rem}\label{rem:12}
Let $\varphi\in \EE_M^{\bar 0}$ be nilpotent.  Then $\ad_\varphi$ is also nilpotent and $\Ad_\varphi=e^{\ad_\varphi}$. Furthermore, since
\begin{equation}\label{eq:8}
\Ad_\varphi \circ \ad_\psi \circ \Ad_{-\varphi} = \ad_{\Ad_\varphi(\psi)}
\end{equation}
for all $\psi\in \EE_M$, we obtain
\begin{equation}\label{eq:9}
\Ad_\varphi(\llbracket \psi_1,\psi_2\rrbracket_\delta)=\llbracket \Ad_\varphi(\psi_1),\Ad_\varphi(\psi_2)\rrbracket_{\Ad_\varphi(\delta)}
\end{equation}
for all $\psi_1,\psi_2,\delta\in \EE_M$.
\end{rem}

\begin{mydef}\label{def:13}
The ${\it Leibnizator}$ of $\delta\in \EE_M$ is the function ${\rm \mathcal L}_\delta:\EE_M\times \EE_M\to \EE_M$ defined by ${\rm \mathcal L}_\delta(\varphi,\psi)=[\L_\varphi^{\!\!^\delta},\L_\psi^{\!\!^\delta}]-\L_{{\llbracket \varphi,\psi\rrbracket}_\delta}^{\!\!^\delta}$ for each $\varphi,\psi\in \EE_M$.
\end{mydef}

\begin{lem}\label{lem:14}
If $\delta\in \EE_M^{\bar 1}$, then
$2{\mathcal L}_\delta(\varphi,\psi)=-(-1)^k\ad_{\llbracket \varphi,\psi\rrbracket_{\delta^2}}$ for all $\varphi\in \EE_M$ and $\psi\in\EE_M^k$. 
\end{lem}

\smallskip\noindent\emph{Proof:} Using \eqref{eq:3} repeatedly, we obtain
\begin{equation}\label{eq:10}
2\mathcal L_\delta(\varphi,\psi)=2\ad_{[[\varphi,\delta],[\psi,\delta]]-[[[\varphi,\delta],\psi],\delta]}=-(-2)^k\ad_{[[[\varphi,\delta],\delta],\psi]}=-(-1)^k\ad_{\llbracket \varphi,\psi\rrbracket_{\delta^2}}\,.
\end{equation}

\begin{example}\label{ex:15}
In particular, $d^2=0$ implies that the Leibnizator of the de Rham operator vanishes identically. Unraveling the definition we obtain 
\begin{equation}\label{eq:11}
\llbracket\varphi_1,\llbracket \varphi_2,\varphi_3 \rrbracket \rrbracket = \llbracket\llbracket \varphi_1,\varphi_2 \rrbracket , \varphi_3\rrbracket+(-1)^{k_1+k_2+k_1k_2}\llbracket \varphi_2,\llbracket \varphi_1,\varphi_3\rrbracket\rrbracket  
\end{equation}
for all $\varphi_i\in \EE_M^{k_i}$, $i=1,2,3$. In particular, if $\phi_1,\phi_2,\phi_3\in \Gamma(\T M)$ we recover the familiar Jacobi identity for the Dorfman bracket.
\end{example}

\section{Generalized Almost Product Structures}

\begin{mydef}\label{def:16}
A {\it split structure of rank $2k$ on $M$} is a subbundle $E\subseteq \T M$ on which the restriction $\langle\,,\,\rangle_E$ of the tautological bilinear form $\langle\, ,\,\rangle$ to $E$ is non-degenerate and of signature $(k,k)$. We denote by $\pr_E$ the orthogonal projection of $\T M$ onto $E$.
\end{mydef}

\begin{rem}\label{rem:17}
If $B\in \Omega_M^2$, then $\Ad_B$ restricts to an automorphism of $\T M$, which by \eqref{eq:8} is orthogonal with respect to $\langle\,,\,\rangle$. In particular, $E\subseteq \T M$ is a split structure if and only if $\Ad_B(E)$ is.
\end{rem}

\begin{rem}\label{rem:18}
$E$ is a split structure of rank $2k$ on $M$ if and only if $E^\perp$ is a split structure of rank $2(\dim(M)-k)$.
\end{rem}

\begin{rem}\label{rem:19}
If $E$ is a non-zero split structure on $M$, non-degeneracy implies that  $\pi(E)$ has nowhere vanishing fibers. Using partitions of unity and the local existence theorem for ODEs, it follows that every function in $\Omega_M^0$ can be written locally as $\L_{\pi(x)}(f)$ for some $x\in \Gamma(E)$ and for some $f\in \Omega_M^0$. On the other hand, \eqref{eq:4} implies
\begin{align}\label{eq:12}
(\ad_\varphi\circ\L_{\pi(x)})(f) & =\ad_\varphi([x,df])=\ad_\varphi([x,\pr_E(df)])\nonumber\\
&= -(-1)^k(\llbracket x,\pr_E(df)\rrbracket_\varphi+\llbracket \pr_E(df),x\rrbracket_\varphi)
\end{align}
for $x\in \Gamma(E)$, $f\in \Omega_M^0$ and for any $\varphi\in \EE_M^k$. Setting $\varphi=d$, we conclude that $E$ is closed under the Dorfman bracket if and only if $T^*M\subseteq E$ if and only if $E=\T M$. 
\end{rem}

\begin{lem}\label{lem:20}
Let $\varphi\in \EE_M$ and assume there exists a split structure $E$ on $M$ such that $\llbracket x,y\rrbracket_\varphi\in \ker (\ad)$ for all $x,y\in \Gamma(E)$. Then $\varphi$ is $\Omega_M^0$-linear.
\end{lem}

\smallskip\noindent\emph{Proof:}
By assumption, $\llbracket x,y\rrbracket_\varphi$ commutes with $\Gamma(E)$ and with $d$ and thus must be a constant multiple of the identity of all $x,y\in \Gamma(E)$. By \eqref{eq:3}, we conclude that $\ad_\varphi$ is a derivation of $\Omega_M^0$ whose image consists of constant functions. This concludes the proof since the only such derivation is the zero derivation.

\begin{lem}\label{lem:21}
Let $\delta\in \EE_M^{\bar 1}$ be such that $\delta(1)=0$ and $[\delta^2,d]=0$. Then following are equivalent
\begin{enumerate}[1)]
\item $\delta^2=0$;
\item $\mathcal L_\delta=0$;
\item there exists a non-zero split structure $E$ on $M$ such $\mathcal L_\delta (x,y)=0$ for all $x,y\in \Gamma(E)$.
\end{enumerate}
\end{lem}

\smallskip\noindent\emph{Proof:}
By Lemma \ref{lem:14}, 1) implies 2). Since $\T M$ is a split structure, 2) implies 3). If 3) holds, then it follows from Lemma \ref{lem:14} and Lemma \ref{lem:20} that $\delta^2$ is $\Omega_M^0$-linear. Combining Remark \ref{rem:8} with the assumption $[\delta^2,d]=0$, we conclude that $\delta^2=\delta^2(1)=0$ which concludes the proof.

\begin{mydef}\label{def:22}
Let $E$ be a split structure of rank $2k$ on $M$. The {\it type of $E$ at $m\in M$} is the rank, denoted by $p_E(m)$, of $\pi_E$ at $m$. 
\end{mydef}

\begin{rem}\label{rem:23} 
Let $E$ be a split structure of rank $2k$ on $M$.
Since $\ker(\pi_E)=T^*M\cap E$ is isotropic, then $p_E(m)\ge k$ for all $m\in M$.
\end{rem}

\begin{mydef}\label{def:24}
A {\it generalized almost product structure} is a split structure $E$ of rank $2k$ on $M$ such that $p_E=k$.
\end{mydef}

\begin{example}\label{ex:25}
Recall that an {\it almost product structure} \cite{gugenheim-spencer} is a pair $(F,G)$ of subbundles of $TM$ such that $TM=F\oplus G$. Each almost product structure $(F,G)$ defines two canonical generalized almost product structures: $E=F\oplus {\rm Ann}(G)$ and $E^\perp=G\oplus {\rm Ann}(F)$.
\end{example}

\begin{prop}\label{prop:26}
Let $E$ be a split structure on $M$. The following are equivalent
\begin{enumerate}[1)]
\item $E$ is a generalized almost product structure;
\item $T^*M=(T^*M\cap E)\oplus (T^*M \cap E^\perp)$;
\item $(\pi(E),\pi(E^\perp))$ is an almost product structure;
\item There exists $B_E\in \Omega^2_M$ such that
\[
\Ad_{B_E}(E)=\pi(E)\oplus {\rm Ann}(\pi(E^\perp))\quad \textrm{ and }\quad \Ad_{B_E}(E^\perp)=\pi(E^\perp)\oplus {\rm Ann}(\pi(E))\,;
\]
\item $[\pr_E(df),\pr_E(dg)]=0$ each $f,g\in \Omega_M^0$.
\end{enumerate}
\end{prop}

\smallskip\noindent\emph{Proof:} Since $E$ is a generalized almost product structure on $M$, then $T^*M\cap E\subseteq E$ is maximal isotropic. On the other hand, ${\rm Ann}(\pi(E))=T^*M\cap E^{\perp}$ and is a subbundle of rank $n-k$. Consequently,  1) implies 2). If 2) holds, then $T^*M\cap E\subseteq E$ and $T^*M\cap E^\perp \subseteq E^\perp$ are maximal isotropic. Therefore $E$ and $E^\perp$ are generalized almost product structures. Since $TM=\pi(E\oplus E^\perp)=\pi(E)+\pi(E^\perp)$, we conclude that 2) implies 1) and 3). Since 
\begin{equation}\label{eq:13}
0=[x,y]=\ad_y(\pi(x))+\ad_x(\pi(y))\,.
\end{equation}
for any $x\in \Gamma(E)$ and $y\in \Gamma(E^\perp)$, if 3) holds  there exists a well defined $B_E\in \Omega_M^2$ such that $(\ad_{\pi(x)}\circ\ad_{\pi(y)}) (B_E) = \ad_y(\pi(x))$ if $x\in \Gamma(E)$ and $y\in \Gamma(E^\perp)$ while $(\ad_{\pi(x)}\circ\ad_{\pi(y)}) (B_E)=0$ if $x,y\in \Gamma(E)$ or $x,y\in \Gamma(E^\perp)$. If $x\in \Gamma(E)$, then $\pi(\Ad_{B_E}(x))=\pi(x)$ and
\begin{align}
\ad_{\pi(y)}(\Ad_{B_E}(x)-\pi(x))&=\ad_y(x-\pi(x)-\ad_{\pi(x)}(B_E))
\nonumber \\ 
&=-\ad_y(\pi(x))-(\ad_{\pi(y)}\circ\ad_{\pi(x)})(B_E)\\
&=0\label{eq:14}
\end{align}
for all $y\in \Gamma(E^\perp)$. Therefore, $\Ad_{B_E}(E)\subseteq \pi(E)\oplus {\rm Ann}(\pi(E^\perp))$. On the other hand
\begin{equation}\label{eq:15}
(\ad_y\circ\Ad_{-B_E})(\pi(x))=\ad_y(\pi(x))+(\ad_{\pi(y)}\circ\ad_{\pi(x)})(B_E)=0
\end{equation}
for all $x\in \Gamma(E)$ and $y\in \Gamma(E^\perp)$. Together with
\begin{equation}\label{eq:16}
\Ad_{-B_E}({\rm Ann}(\pi(E^\perp))={\rm Ann}(\pi(E^\perp))\subseteq E
\end{equation}
this implies $\pi(E)\oplus {\rm Ann}(\pi(E^\perp))\subseteq \Ad_{B_E}(E)$. A similar calculation for $E^\perp$ shows that 3) implies 4). If 4) holds, then $(\Ad_{B_E}(E))^\perp=\Ad_{B_E}(E^\perp)$ implies
\begin{equation}\label{eq:17}
T^*M=(T^*M\cap \Ad_{B_E}(E))\oplus (T^*M\cap \Ad_{B_E}(E^\perp))=(T^*M\cap E)\oplus (T^*M \cap E^\perp)
\end{equation}
and thus 1). Finally for each fixed $f\in \Omega_M^0$, $[\pr_E(df),dg]=[\pr_E(df),\pr_E(dg)]=0$
for all $g\in \Omega_M^0$ if and only if $\pr_E(df)\in \Gamma(T^*M\cap E)$. Therefore, 5) is equivalent to 2). 

\begin{cor}\label{cor:27}
A split structure $E$ on $M$ is a generalized almost product structure if and only if $E^\perp$ is.
\end{cor}

\begin{mydef}\label{def:28}
Let $E$ be an almost product structure on $M$. The {\it $E$-bigrading} is the canonical $(\Z\times \Z)$-grading of $\Omega_M$ induced by the decomposition $T^*M=(T^*M\cap E)\oplus (T^*M\cap E^\perp)$ so that the bigraded components are $\Omega_E^{k,l}=\Gamma(\wedge^k(T^*M\cap E)\otimes\wedge^l(T^*M\cap E^\perp))$. We denote by $\EE_E^{k,l}$ the components of the induced decomposition of $\EE_M$. The {\it $E$-biparity} is the $(\Z/2\times\Z/2)$-reduction of the $E$-bigrading. The {\it $E$-parity} is the $\Z/2$-grading corresponding to the decomposition $\Omega_M=\Omega_E^{\overline 0,\bullet}\oplus \Omega_E^{\overline 1,\bullet}$.
\end{mydef}

\begin{rem}\label{rem:29}
The $E$-biparity is a simultaneous $(\Z/2\times\Z/2)$-refinement of both the standard parity and the $E$-parity.
\end{rem}

\begin{cor}\label{cor:30}
Generalized almost product structures on $M$ are in canonical bijection with pairs $(\Omega^{\bullet,\bullet},B)$, where $\Omega^{\bullet,\bullet}$ denotes a $(\Z\times \Z)$-refinement of the standard grading on $\Omega_M$ and $B\in \Omega^{1,1}$.
\end{cor}

\begin{rem}\label{rem:31}
Let $E$ be a generalized almost product structure on $M$. If $B_E\neq 0$, then sections of $E$ (and of $E^\perp)$ do not in general have definite $E$-bigrading. However, Proposition \ref{prop:26} implies that $\Gamma(E)\subseteq \EE_E^{1,0}\oplus \EE_E^{-1,0}\oplus \EE_E^{0,1}$.
\end{rem}

\begin{rem}\label{rem:32}
Let $E$ be a generalized almost product structure on $M$. By construction, the $E$-bigrading coincides with the canonical bigrading of the almost product structure $(\pi(E),\pi(E^\perp))$, as defined in \cite{gugenheim-spencer}. In particular, $d\in \EE_E^{1,0}\oplus \EE_E^{-1,2}\oplus \EE_E^{0,1}\oplus \EE_E^{2,-1}$.
\end{rem}

\begin{mydef}\label{def:33}
Given a generalized almost product structure $E$ on $M$, we define $d_E$ to be the component of odd $E$-parity of the de Rham operator, so that $d_E\in \EE_E^{1,0}\oplus \EE_E^{-1,2}$. We use the shorthand notations $\L^E=\L^{d_E}$ and $\llbracket\, ,\, \rrbracket_E=\llbracket\, ,\, \rrbracket_{d_E}$.
\end{mydef}

\begin{example}\label{ex:34}
Let $M=S^3$, let $\{X_1,X_2,X_3\}$ be a frame of $TM$ with dual frame $\{\alpha_1,\alpha_2,\alpha_3\}$ such that $d\alpha_1=\alpha_2\alpha_3$, $d\alpha_2=\alpha_3\alpha_1$ and $d\alpha_3=\alpha_1\alpha_2$. Given $f,g\in \Omega_M^0$, consider the split structure $E={\rm span}\{\alpha_2,\alpha_3,x_2,x_3\}$ where $x_2=X_2-f\alpha_1$ and $x_3=X_3-g\alpha_1$. It follows that $E^\perp={\rm span}\{\alpha_1,x_1\}$, where $x_1=X_1+f\alpha_2+g\alpha_3$. Then $E$ is a generalized almost product structure for all $f,g\in \Omega_M^0$. A direct calculation shows that $d_E(h\alpha_i)=d_E(h)\alpha_i$ for $i=1,2,3$.
\end{example}

\begin{rem}\label{rem:35}
If $E$ is a generalized almost product structure on $M$, then $d_Ef=\pr_E(df)$ for every $f\in \Omega_M^0$. Moreover, $d_E=d_{\Ad_B(E)}$ for every $B\in \Omega_M^2$.
\end{rem}

\begin{rem}\label{rem:36} 
Let $E$ be a generalized almost product structure on $M$. Separating the terms of different $E$-biparity in the identity $d^2=0$ yields $[d_E,d_{E^\perp}]=0$ and $d_E^2+d_{E^\perp}^2=0$. In particular, $d_E^2=0$ if and only if $d_{E^\perp}^2=0$ if and only if $[d,d_E]=0$. On the other hand
\begin{equation}\label{eq:18}
[d_E,[d_E,d_E]]=-[d_E,[d_{E^\perp},d_{E^\perp}]]=-2[d_{E^\perp},[d_E,d_{E^\perp}]]=0
\end{equation}
for any generalized almost product structure $E$.
\end{rem}

\begin{rem}\label{rem:37}
Let $E$ be a generalized almost product structure on $M$ and let $x,y\in \Gamma(E)$. Then $\ad_{\llbracket x,y\rrbracket_E}=[[\ad_x,\ad_{d_E}],\ad_y]$ is an operator of $E$-biparity $(\bar 1,\bar 0)$ on $\EE_M$. Similarly, $\ad_{\llbracket x,y\rrbracket_{E^\perp}}$ has $E$-biparity $(\bar 0,\bar 1)$. Matching $E$-biparities yields $\llbracket x,y\rrbracket_E=\pr_E(\llbracket x,y\rrbracket)$ and $\llbracket x,y\rrbracket_{E^\perp}=\pr_{E^\perp}(\llbracket x,y\rrbracket)$. Similarly if $x\in \Gamma(E)$ and $y\in \Gamma(E^\perp)$, then $\llbracket x,y\rrbracket_E=\pr_{E^\perp} (\llbracket x,y\rrbracket)$.   
\end{rem}

\begin{rem}\label{rem:38}
Let $E$ be a generalized almost product structure on $M$ and let $x,y\in \Gamma(E)$. Since by Lemma \ref{lem:14} the components of different $E$-biparity in $\mathcal L_d(x,y)$ must vanish independently, it follows that in particular $\mathcal L_E(x,y)=-\mathcal L_{E^\perp}(x,y)$ for every $x,y\in \Gamma(E)$.
\end{rem}

\begin{rem}\label{rem:39}
Let $E$ be a generalized almost product structure on $M$. By projecting onto $E$ the axioms of Courant algebroid (written in terms of the Dorfman bracket as in \cite{G}), we conclude that $(E,\llbracket\,,\,\rrbracket_{_E},\langle\,,\,\rangle_{_E},\pi_E)$ is a Courant algebroid if and only if $\llbracket\,,\,\rrbracket_{_E}$ satisfies the Jacobi identity i.e. $\Gamma(E)\subseteq \ker (\mathcal L_E(x,y))$ for all $x,y,\in \Gamma(E)$.
\end{rem}

\begin{rem}\label{rem:40}
Let $E$ be a generalized almost product structure on $M$. By a result of \cite{reinhart}, $\pi(E)$ is a foliation if and only if $d_{E^\perp}\in \EE_E^{0,1}$. If this is the case, we see that in particular $d_E^2=-d_{E^\perp}^2\in \EE_E^{0,2}$.
\end{rem}

\begin{prop}\label{prop:41}
Let $E$ be an almost product structure on $M$. The following are equivalent
\begin{enumerate}[1)]
\item $\pi(E)$ is a foliation;
\item $\Omega_M^0\subseteq \ker(\mathcal L_E(x,y))$ for all $x,y\in \Gamma(E)$;
\item $\llbracket x,y\rrbracket_{E^\perp}\in \Omega_M^1$ for all $x,y\in \Gamma(E)$;
\item $(E,\llbracket\,,\,\rrbracket_E,\langle\,,\rangle_E,\pi_E)$ is a Courant algebroid.
\end{enumerate}
\end{prop}

\smallskip\noindent\emph{Proof:} By Remark \ref{rem:37}, 3) is equivalent to 
\begin{equation}\label{eq:19}
0=[\llbracket x,y\rrbracket_{E^\perp},df]=[\llbracket x,y\rrbracket_{E^\perp},d_{E^\perp} f]
\end{equation}
for all $x,y\in \Gamma(E)$ and for all $f\in \Omega_M^0$. Since $\llbracket x,f\rrbracket_{E\perp}=[x,d_{E^\perp}f]=0$ for all $x\in \Gamma(E)$ and $f\in \Omega_M^0$, using Remark \ref{rem:38} we obtain
\begin{equation}\label{eq:20}
(\mathcal L_E(x,y))(f)=-(\mathcal L_{E^\perp}(x,y))(f)=\L^{\!\!^{E^\perp}}_{\llbracket x,y\rrbracket_{E^\perp}}(f)=[\llbracket x,y\rrbracket_{E^\perp},d_{E^\perp} f]\,.
\end{equation}
Therefore, 2) is equivalent to 3). On the other hand, Lemma \ref{lem:14} implies
\begin{equation}\label{eq:21}
2(\mathcal L_E(x,y))(f)=[\llbracket x,y\rrbracket_{d_E^2},f]=\llbracket x,y\rrbracket_{d_E^2f}
\end{equation}
for all $x,y\in \Gamma(E)$ and for all $f\in \Omega_M^0$. By Lemma \ref{lem:14}, 4) is equivalent to 
\begin{equation}\label{eq:22}
0=[\llbracket x,y\rrbracket_{d_E^2},z]=\llbracket y,z\rrbracket_{[x,d_E^2]}
\end{equation}
for all $x,y,z\in \Gamma(E)$. Lemma \ref{lem:20} then implies that $[x,d_E^2]$ is $\Omega_M^0$-linear and (since $x$ and $y$ are also $\Omega_M^0$-linear) we obtain that $\llbracket x,y\rrbracket_{d_E^2f}=0$ for all $x,y\in \Gamma(E)$ and for all $f\in \Omega_M^0$. Therefore, 4) implies 2). Since 4) clearly implies that $\pi(E)$ is a foliation, it remains to prove that 1) implies 4). Let $x,y,z$ be arbitrary sections of $E$. On the one hand, $(\mathcal L_E(x,y))(z)$ has only components of bidegree $(\bullet,q)$ with $q\le 1$ by Remark \ref{rem:31}. On the other hand by by Remark \ref{rem:38} and Remark \ref{rem:40}, $(\mathcal L_E(x,y))(z)$ has only components of $E$-bigradee $(\bullet,q)$ with $q\ge 2$. Therefore, it must vanish and the proof is completed.

\begin{mydef}\label{def:42}
Let $E$ be a generalized almost product structure on $M$ such that $\pi(E^\perp)$ is a foliation. A differential form is {\it basic with respect to $E$} if it is an element of $\mathcal B_E = \Omega^{\bullet,0}_E\cap \ker(d_{E^\perp})$. Accordingly, an operator $\varphi\in \EE_M$ is {\it basic with respect to $E$} if $\varphi(\mathcal B_E)\subseteq \mathcal B_E$.
\end{mydef}

\begin{example}\label{ex:43}
Let $E$ be a generalized almost product structure on $M$ such that $\pi(E^\perp)$ is a foliation, then $d_E$ is basic with respect to $E$. Moreover $d_E^2(\mathcal B_E)=d_{E^\perp}^2(\mathcal B_E)=0$ and thus $(\mathcal B_E,d_E)$ is a complex known as the {\it basic complex of $E$}. 
\end{example}

\begin{mydef}\label{def:44}
A split structure $E$ is a {\it generalized local product structure} on $M$ if $(\pi(E),\pi(E^\perp))$ is a local product structure in the sense of \cite{reinhart} i.e.\ $TM=\pi(E)\oplus \pi(E^\perp)$ is a decomposition into constant-rank foliations.
\end{mydef}

\begin{prop}\label{prop:45}
A split structure $E$ on $M$ is a generalized local product structure if and only if it is a generalized almost product structure and $d_E^2=0$. 
\end{prop}

\smallskip\noindent\emph{Proof:} By Proposition \ref{prop:26}, every generalized local product structure is also a generalized almost product structure. By Remark \ref{rem:40}, if $(\pi(E),\pi(E^\perp))$ is a local product structure, then $d_E^2\in \EE_E^{0,2}\cap \EE_{E^\perp}^{0,2}=\EE_E^{0,2}\cap \EE_E^{2,0}=\{0\}$. Conversely, if $d^2_E=0$ then by \eqref{eq:22}, $\pi(E)$ is a foliation. By Remark \ref{rem:40}, $\pi(E^\perp)$ is also a foliation and the Proposition is proved. 

\begin{cor}\label{cor:46}
The collection of all split structures on $M$ such that $d^2_E=0$ is in canonical bijection with the collection of pairs $(\Omega^{\bullet,\bullet},B)$ consisting of a $(\Z\times \Z)$-refinement of the standard grading on $\Omega_M$ such that $d(\Omega^{i,j})\subseteq \Omega^{i+1,j}\oplus \Omega^{i,j+1}$ for all $i,j\ge 0$ and $B\in \Omega^{1,1}$.
\end{cor}

\section{Generalized $F$-structures}

\begin{mydef}\label{def:47}
Let $E$ be a split structure on $M$. A {\it generalized $F$-structure on $E$} is an element $\Phi\in \EE_M$ such that
\begin{enumerate}[i)]
\item $\Omega_M^0\oplus \Gamma(E^\perp)\subseteq \ker(\ad_\Phi)$;
\item $\ad_\Phi(\Gamma(\T M))\subseteq \Gamma(\T M)$;
\item $\Gamma(E)\subseteq \ker(\ad_\Phi^2+{\rm Id})$.
\end{enumerate}
We denote by $J_\Phi$ the restriction of $\ad_\Phi$ to $\Gamma(\T M)$ and by $L_\Phi$ the $\sqrt{-1}$-eigenbundle of $J_\Phi$.
\end{mydef}

\begin{rem}\label{rem:48}
If $\Phi$ is a generalized $F$-structure on a split structure $E$, $J_\Phi$ is an orthogonal bundle endomorphism of $\T M$ such that $J_\Phi^3+J_\Phi=0$ i.e.\ a generalized $F$-structure in the sense of \cite{vaismanCRF}. Conversely, given any bundle endomorphism $J\in {\mathfrak o} (\T M)$, there is an element $\Phi\in \EE_M$ such that $J(x)=\ad_\Phi(x)$ for all $x\in \T M$. For instance, let $\omega$ be the unique $2$-form on $M$ such that $\omega(X,Y)=[X,J(Y)]$ for all $X,Y\in \Gamma(TM)$ and consider $\Phi\in \EE_M$ such that $\Phi(f)=f\omega$ for any $f\in \Omega_M^0$ and 
\begin{equation}\label{eq:23}
\Phi(\alpha_1\cdots\alpha_p)=\alpha_1\cdots\alpha_k\omega+\sum_{i=1}^p \alpha_1\cdots J(\alpha_i)\cdots \alpha_k
\end{equation}
for any $\alpha_1,\ldots,\alpha_p\in \Omega^1_M$. Then $\Phi$ is $\Omega_M^0$-linear and the restriction of $\ad_\Phi$ to $\T M$ coincides with $J$. Furthermore, suppose that $\Phi'\in \EE_M$ is another $\Omega_M^0$-linear element such that $J(x)=\ad_{\Phi'}(x)$, then Remark \ref{rem:3} implies $\Phi-\Phi'\in\Omega_M^0$. Therefore, modulo addition of functions, generalized $F$-structures on split structures	 are in canonical correspondence with the split generalized $F$-structures defined in \cite{AG2}.
\end{rem}

\begin{example}\label{ex:49}
Every generalized almost complex structure is of the form $J_\Phi$ for some generalized $F$-structures $\Phi$ on $M$.
\end{example}

\begin{example}\label{ex:50}
In the language of \cite{AG}, every generalized almost contact triple is of the form $(J_\Phi,e_1,e_2)$ for some generalized $F$-structure $\Phi$ on a split structure $E$ such that $E^\perp$ is globally trivialized by isotropic sections $e_1,e_2\in \Gamma(E^\perp)$. If one further imposes the condition $e_1\in \Gamma(TM),e_2\in \Gamma(T^*M)$ one obtains the generalized almost contact triples of
\cite{PW}.
\end{example}

\begin{rem}\label{rem:51}
Since $J_\Phi^3+J_\Phi=0$, then $L_\Phi$ is empty if and only if $\Phi$ is $\Omega_M^0$-linear. Moreover, $L_\Phi\subseteq E\otimes \mathbb C$ is maximal isotropic with respect to the tautological inner product and $E\otimes \C=L_\Phi\oplus\overline L_\Phi$.  
\end{rem}

\begin{example}\label{ex:52}
Let $E$ be a generalized almost product structure on $M$, let $\omega\in \Omega_E^{2,0}$ and let $\Lambda\in \EE_E^{-2,0}$ be such that $\Phi=\omega+\Lambda$ is a generalized $F$-structure on $E$. Then
\begin{equation}\label{eq:24}
\Gamma(L_\Phi)=({\rm Id}-\sqrt{-1}\ad_\Lambda)(\Omega^{1,0}_E)\,.
\end{equation}
\end{example}

\begin{lem}\label{lem:53}
Let $E$ be a split structure on $M$ and let $\Phi$ be a generalized $F$-structure on $E$. The following are equivalent:
\begin{enumerate}[1)]
\item $J_\Phi(E)\subseteq E$;
\item $L_\Phi = (L_\Phi\cap (E\otimes \C))\oplus (L_\Phi\cap (E^\perp\otimes \C))$;
\item there exists a generalized $F$-structure $\Phi_E$ on $E$ such that $\Phi-\Phi_E$ is a generalized $F$-structure on $E^\perp$.
\end{enumerate}
\end{lem}

\smallskip\noindent\emph{Proof:} The implications $3)\Rightarrow 1)\Rightarrow 2)$ are straightforward. To see why 2) implies 3), let $\{l_1,\ldots,l_n\}$ be a local frame of $L_\Phi$ such that $l_1,\ldots,l_k\in E\otimes \C$ and $l_{k+1},\ldots,l_n\in E^\perp\otimes \C$. If $\{\overline l^1,\ldots,\overline l^n\}\subseteq \Gamma(\overline L_\Phi)$ is the dual local frame, then $\sqrt{-1}\sum_{i=1}^n (l_i\circ \overline l^i)$ is a local expression for $\Phi$ and $\sqrt{-1}\sum_{i=1}^k (l_i\circ \overline l^i)$ is a local expression for $\Phi_E$. 

\begin{mydef}\label{def:54}
Let $E$ be a generalized almost product structure on $M$ and let $\Phi$ be a generalized $F$-structure on $E$. Let $K'_\Phi\subseteq \wedge^\bullet T^*M\otimes \C$ be the pure spinor line bundle of the maximal isotropic subbundle $L_\Phi\oplus \Ad_{-B_E}(\pi(E^\perp))\otimes \C\subseteq \T M\otimes \C$. The {\it canonical bundle of $\Phi$} is $K_\Phi=K'_\Phi\otimes\wedge^\bullet (T^*M\cap E^\perp)\subseteq \wedge^\bullet T^*\otimes \C$. 
\end{mydef}

\begin{rem}\label{rem:55}
Let $\Phi$ be a generalized $F$-structure on a generalized almost product structure $E$ of rank $2k$ on $M$. Then $K_\Phi$ is a complex bundle of  rank $2^{\dim M-k}$.
\end{rem}

\begin{rem}\label{rem:56}
Let $E$ be a generalized almost product structure on $M$ and let $\Phi$ be a generalized $F$-structure on $E$.
Since $L_\Phi\oplus\Ad_{-B_E}(\pi(E^\perp))\otimes \C ={\rm Ann}(K'_\Phi)$ and $L_\Phi\subseteq (T^*M\cap E^\perp)\otimes \C$, then $L_\Phi\subseteq {\rm Ann}(K_\Phi)\subseteq L_\Phi\oplus \Ad_{-B_E}(\pi(E^\perp))\otimes \C$. On the other hand, no section of $\Ad_{-B_E}(\pi(E^\perp))$ annihilates every section of $T^*M\cap E^\perp$ and thus $L_\Phi={\rm Ann}(K_\Phi)$.
\end{rem}

\begin{example}\label{ex:57}
If $\Phi$ is a generalized $F$-structure on $\T M$, then $K_\Phi=K'_\Phi$ is the canonical line bundle of the generalized almost complex structure $J_\Phi$, in the sense of \cite{G}.
\end{example}

\begin{rem}\label{rem:58}
Let $E$ be a generalized almost product structure of rank $2k$ on $M$ and let $\Phi$ be a generalized $F$-structure on $E$. For each integer $i\ge 0$, let $F_\Phi^i\subseteq \wedge^\bullet T^*M\otimes \C$ be the subbundle annihilated by $\Gamma(\wedge^{i+1} L_\Phi)$. In particular, $F_\Phi^0=K_\Phi$ and $F_\Phi^i=\wedge ^\bullet T^*M\otimes \C$ for all $i\ge 2j$. Since $\Gamma(F^{i+1}_\Phi)\setminus \Gamma(F_\Phi^i)$ contains all images of forms in $\Gamma(K_\Phi)$ under the actions of operators in $\Gamma(\wedge^{i+1} \overline L_\Phi)$, then $F^{i+1}_\Phi/F^i_\Phi$ is a bundle of rank at least $\binom{2k}{i}$ which yields a canonical isomorphism
\begin{equation}\label{eq:25}
F^i_\Phi\cong \bigoplus_{j=0}^i (\wedge^j \overline L_\Phi)\otimes K_\Phi\,.
\end{equation}
Since $\overline F_\Phi^i=F^i_{-\Phi}$ one obtains a further canonical isomorphism
\begin{equation}\label{eq:26}
F_\Phi^i\cap \overline F_\Phi^{2k-i}\cong \wedge^i \overline L_\Phi \otimes K_\Phi\,.
\end{equation}
Letting $U^{k-i}_\Phi=\Gamma(F^i_\Phi\cap \overline F^{2k-i}_\Phi)$ we obtain a new $\Z$-grading, called the {\it $\Phi$-grading},
\begin{equation}\label{eq:27}
\Omega_M\otimes \C=U^{-k}_\Phi\oplus\cdots \oplus U^k_\Phi\,.  
\end{equation}
Note that $U_\Phi^k=\Gamma(K_\Phi)$ and $\overline U_\Phi^i=U_\Phi^{-i}$ for all $i$. The component of degree $i$ of the corresponding $\Z$-grading of $\EE_M\otimes \C$ is denoted by $\EE_\Phi^i$.
\end{rem}

\begin{example}\label{ex:59}
Let $\Phi$ be a generalized $F$-structure on $\T M$. Then the $\Phi$-grading defines the standard decomposition of complex differential forms associated to the generalized almost complex structure $J_\Phi$, as defined in \cite{G}. 
\end{example}

\begin{example}\label{ex:60}
Let $\Phi=\omega+\Lambda$ be as in Example \ref{eq:52} and assume that $E$ is of rank $2k$. Then $U_\Phi^k=e^{-\sqrt{-1}\omega}\Omega_{E}^{0,\bullet}$. Moreover for every $\alpha\in \Omega^{1,0}_E$
\begin{equation}\label{eq:28}
\Psi\circ \alpha\circ \Psi^{-1}=\left(\Ad_{\sqrt{-1}\omega}\circ\Ad_{-\frac{\sqrt{-1}}{2}\Lambda}\right)(\alpha)=-\frac{1}{2}(\alpha+\sqrt{-1}\ad_\Lambda(\alpha))\,,
\end{equation}
where $\Psi=e^{\sqrt{-1}\omega}\circ e^{-\frac{\sqrt{-1}}{2}\Lambda}$. Therefore, $\Psi$ intertwines the standard action of forms with the action of $\overline L_\Phi$ and thus gives rise to a canonical identification $U^{k-i}_\Phi=\Psi(\Omega_E^{i,\bullet})$ for all $i$. In the particular case in which $\omega$ is non-degenerate and $J_\Phi$ is the corresponding generalized almost complex structure, we obtain Theorem 2.2 in \cite{cavalcanti}.
\end{example}

\begin{rem}\label{rem:61}
Specializing Example \ref{ex:60} to the case in which $\dim M=2k+1$ and $(\omega,\eta)$ is an almost cosymplectic structure i.e.\ $\eta\in \Omega^1_M$ is such that $\eta\omega^k$ is a volume form, then $U_\Phi^{-k}$ is the $\Omega_M^0$-module generated by $e^{\sqrt{-1}\omega}$ and $\eta e^{\sqrt{-1}\omega}$. In particular we observe that, unlike the special case of generalized almost complex structures, in general the subspaces $U^i_\Phi$ do not have definite standard parity. 
\end{rem}

\begin{rem}\label{rem:62}
Let $E$ be a generalized almost product structure of rank $2k$ on $M$ and let $\Phi$ be a generalized $F$-structure on $E$. By Remark \ref{rem:58}, $U^k_\Phi$ is a subspace of $\Omega_E^{\bullet,0}$ generated by pure spinors of definite standard parity. Since the standard parity and the $E$-parity coincide when restricted to $\Omega_E^{\bullet,0}$, we conclude that elements of $U^{k}_\Phi$, and thus of $U^i_\Phi$ for all $i$, have definite $E$-parity. More precisely, the $\Z/2$-reduction of the $\Phi$-grading coincides with the $E$-parity.
\end{rem}

\begin{lem}\label{lem:63}Let $E$ be a generalized almost product structure on $M$, let $\Phi$ be a generalized $F$-structure on $E$.
and let $\varphi\in \EE_M$. Then $\varphi\in \EE_\Phi^i$ if and only if $\ad_\Phi(\varphi)=\sqrt{-1}i\varphi$.
\end{lem}

\smallskip\noindent\emph{Proof:} Since $\Phi$ preserves $K'_\Phi$ and commutes with the action of $\Gamma(E^\perp)$ on $\Omega_M$, locally, there exists $f\in \Omega_M^0\otimes \C$ such that $\Phi(\rho)=f\rho$ for each locally defined section $\rho$ of $U_\Phi^k$. Taking into account that $\overline L_\Phi$ is the $-\sqrt{-1}$-eigenbundle of $\ad_\Phi$, we conclude that locally $U_\Phi^{k-i}$ is the $f-i\sqrt{-1}$ eigenspace of $\Phi$ from which the Lemma easily follows.

\begin{lem}\label{lem:64}
Let $E$ be a generalized almost product structure on $M$ and let $\Phi$ be a generalized $F$-structure on $E$. Then 
\begin{enumerate}[1)]
\item $d_{E^\perp}\in \EE_\Phi^{-2}\oplus \EE_\Phi^0\oplus \EE_\Phi^2$;
\item $d_{E^\perp} \in \EE_\Phi^0$ if and only if $\llbracket l_1,l_2\rrbracket_{E^\perp}=0$ for all $l_1,l_2\in \Gamma(L_\Phi)$.
\end{enumerate}
\end{lem}

\smallskip\noindent\emph{Proof:} Assume $E$ has rank $2k$. Since sections of $E^\perp$ commute with $\Phi$, then $\llbracket l_1,l_2\rrbracket_{E^\perp}(U_\Phi^i)\subseteq U_\Phi^i$ for all $l_1,l_2\in \Gamma(L_\Phi)$ and for any integer $i$. Unraveling the $i=k$ case we obtain $d_{E^\perp}(U_\Phi^k)\subseteq \Gamma(F_\Phi^2)$ which, upon inspection of $E$-biparities, yields $d_{E^\perp}(U_\Phi^k)\subseteq U_\Phi^k\oplus U_\Phi^{k-2}$. Arguing by induction on $k-i$, a similar argument shows that $d_{E^\perp}(U_\Phi^i)\subseteq U_\Phi^{i-2}\oplus U_\Phi^i\oplus U_\Phi^{i+2}$ for all $i$. This implies 1). Let $\partial_\Phi^\perp$ be the projection of $d_{E^\perp}$ onto $\EE_\Phi^2$. Matching $\Phi$-degrees, we obtain $\llbracket l_1,l_2\rrbracket_{E^\perp}=\llbracket l_1,l_2\rrbracket_{\overline\partial_\Phi^\perp}$ for all $l_1,l_2\in \Gamma(L_\Phi)$. An immediate consequence is that $d_{E^\perp}\in \EE_\Phi^0$ implies $\llbracket l_1,l_2\rrbracket_{E^\perp}=0$ for all $l_1,l_2\in \Gamma(L_\Phi)$. Conversely, if the restriction of $\llbracket \, ,\, \rrbracket_{E^\perp}$ to $\Gamma(L_\Phi)$ vanishes, Then in particular
\begin{equation}\label{eq:29}
(l_1\circ l_2 \circ \overline\partial_\Phi^\perp)(U^k_\Phi)=\llbracket l_1,l_2\rrbracket_{E^\perp}(U_\Phi^k)=0
\end{equation}
for all $l_1,l_2\in \Gamma(L_\Phi)$ and thus $\overline\partial_\Phi^\perp(U_\Phi^k)\subseteq U_\Phi^k\oplus U_\Phi^{k-1}$. Since $\overline\partial_\Phi^\perp\in \EE_\Phi^{-2}$, we conclude that $\overline\partial_\Phi^\perp U_\Phi^k=0$. Using induction on $k-i$, a similar argument shows $\overline\partial_\Phi^\perp(U_\Phi^i)=0$ for all $i$ i.e. $\partial_\Phi^\perp=\overline\partial_\Phi^\perp=0$. Since $d_{E^\perp}$ is a real operator, it must then be of $\Phi$-degree $0$.

\begin{rem}\label{rem:65}Let $E$ be a generalized almost product structure $E$ on $M$ such that $\pi(E^\perp)$ is a foliation and let $\Phi$ be a generalized $F$-structure on $E$. Since $d_{E^\perp}\Phi(\mathcal B_E)\subseteq \Omega^{\bullet,1}$, we obtain that $\Phi$ is basic with respect to $E$ if and only if $\llbracket \Gamma(\pi(E^\perp)),\Phi\rrbracket_{E^\perp}(\mathcal B_E)=0$. On the other hand, if $X\in \Gamma(\pi(E^\perp))$ then $[\llbracket X,\Phi\rrbracket_{E^\perp},\alpha]=0$ for every $\alpha\in \Omega_E^{0,1}$. Since $\Omega_M$ is generated by $\mathcal B_E\otimes \Omega_E^{0,\bullet}$ we conclude that $\Phi$ is basic with respect to $E$ if and only if $\llbracket \Gamma(\pi(E^\perp)),\Phi\rrbracket_{E^\perp}=0$. Using \eqref{eq:4} and the fact that sections of $L_\Phi$ are of the form $x-\sqrt{-1}[\Phi,x]$ for some $x\in \Gamma(E)$,  we obtain that $\Phi$ is basic with respect to $E$ if and only if $\llbracket \Gamma(\pi(E^\perp)),\Gamma(L_\Phi)\rrbracket_{E^\perp}\subseteq \Gamma(L_\Phi)$, if and only if $\llbracket \Gamma(\pi(E^\perp)),\Gamma(L_\Phi\oplus \pi(E^\perp))\rrbracket\subseteq \Gamma(L_\Phi\oplus \pi(E^\perp))$. In particular, if $K'_\Phi$ is locally generated by spinors that are basic with respect to $E$, then $\Phi$ is itself basic with respect to $E$. Conversely, if $\Phi$ is basic with respect to $E$, then the $\Phi$-grading restricts to a grading of basic forms $\mathcal B_E=\mathcal B_\Phi^{-k}\oplus \cdots \oplus \mathcal B_\Phi^k$, where $\mathcal B_\Phi^i=\mathcal B_E \cap \mathcal U_\Phi^i$ for all $i$.
\end{rem}

\begin{example}\label{ex:66}
Let $\Phi=\omega+\Lambda$ be a generalized $F$-structure as in Example \ref{ex:60}. Assume that $\pi(E^\perp)$ is a foliation and $d_{E^\perp}\omega=0$. Then $K'_\Phi$ is globally trivialized by the form $e^{-\sqrt{-1}\omega}$ which is clearly basic with respect to $E$. Therefore, $\Phi$ is basic with respect to $E$ and thus so is $\Lambda$. It follows from \eqref{eq:28} that $[d_{E^\perp},\Psi]=\frac{\sqrt{-1}}{2}(\Psi\circ \ad_\Lambda)(d_{E^\perp})$ and thus $\mathcal B_\Phi^{k-i}=\Psi (\mathcal B_E^i)$. This generalizes an observation made in \cite{wade10} for the special case in which $\omega=d\eta$ for some contact form $\eta$.  
\end{example}

\begin{rem}\label{rem:67}
While the notion of generalized $F$-structure on a split structure is invariant under T-duality, the notion of generalized almost product structure is not. This suggest to generalize the construction of the canonical bundle to split structures $E$ together with a decomposition $E^\perp=D_1\oplus D_2$ into isotropic subbundles. Then $K'_\Phi$ can be taken to be the spinor line annihilated by $\Gamma(L_\Phi\oplus D_1)$ and $K_\Phi=K'_\Phi\otimes \wedge^\bullet D_2$ so that $L_\Phi={\rm Ann}(K_\Phi)$. The construction of the $\Phi$-grading given in Remark \ref{rem:58} can be extended verbatim to this more general setup. 
\end{rem}

\begin{example}\label{rem:68}
Let $(J_\Phi,e_1,e_2)$ be a generalized almost contact triple as in Example \ref{ex:50} and let $D_i$ be the trivial line bundle generated by $e_i$ for $i=1,2$. Then $K'_\Phi$ is locally generated by a spinor $\rho_1$, which together with $\rho_2=e_2\rho_1$ locally generates $K_\Phi$. In the language of \cite{AG}, $(\rho_1,\rho_2)$ is a local mixed pair.
\end{example}

\section{(Weak) generalized CRF-structures}

\begin{mydef}\label{def:69}
A {\it weak generalized CRF-structure} on a generalized almost product structure $E$ is a generalized $F$-structure $\Phi$ on $E$ whose $\sqrt{-1}$-eigenbundle is closed under $\llbracket \,,\,\rrbracket_E$ i.e.\ $\llbracket l_2,l_2\rrbracket_E\in \Gamma(L_\Phi)$ for each $l_1,l_2\in \Gamma(L_\Phi)$.
\end{mydef}

\begin{example}\label{ex:70}
Let $\Phi$ be a generalized $F$-structure such that $J_\Phi$ is a generalized almost complex structure on $M$. Then $\Phi$ is a weak generalized CRF-structure if and only if $J_\Phi$ is a generalized complex structure. In particular, complex and symplectic structures are particular cases of weak generalized CRF-structures. 
\end{example}

\begin{example}\label{ex:71}
Let $E$ be a split structure of rank $2$ globally trivialized by isotropic sections $e_1,e_2\in \Gamma(E)$. Let $\Phi$ be a generalized $F$-structure on $E^\perp$ such that $(J_\Phi,e_1,e_2)$ is a generalized almost contact triple in the sense of \cite{AG}. For $i=1,2$, let $\mathbb C e_i\subseteq \T M$ be complex line bundle generated by $e_i$. By definition \cite{AG}, the triple $(J_\Phi,e_1,e_2)$ is {\it integrable} if there exist $i\in\{1,2\}$ such that $L_\Phi\oplus \C e_i$ is closed under the Dorfman bracket. Assume that $E$ is an almost product structure i.e.\ either $e_1$ or $e_2$ is a $1$-form. Projecting the Dorfman bracket onto $E^\perp$, it is easy to see that the the integrability of $(J_\Phi,e_1,e_2)$ implies that $\Phi$ is a weak generalized CRF structure. In particular, contact, cosymplectic and normal almost contact structures are examples of weak generalized CRF-structures. 
\end{example}

\begin{example}\label{ex:72}
Let $E$ be a generalized almost product structure such that $\pi(E^\perp)$ is a foliation and $B_E=0$. Let $\Phi$ a generalized $F$-structure on $E$ that is basic with respect to $E$. According to \cite{wade10} $J_\Phi$ is a {\it transverse generalized complex structure} if $\llbracket \Gamma(L_\Phi),\Gamma(L_\Phi)\rrbracket \subseteq \Gamma(L_\Phi\oplus \pi(E^\perp))$. Clearly, this condition implies that $\Phi$ is a weak generalized CRF-structure. Conversely, since $\pi(E^\perp)$ is maximal isotropic in $E^\perp$, the condition $\llbracket \Gamma(L_\Phi),\Gamma(L_\Phi)\rrbracket_{E^\perp}\subseteq\Gamma(\pi(E^\perp))$ is equivalent to $[\llbracket \Gamma(L_\Phi),\Gamma(L_\Phi)\rrbracket_{E^\perp},\Gamma(\pi(E^\perp))]=0$. Using \eqref{eq:4} and the maximal isotropy of $L_\Phi$ in $E\otimes \C$ this is in turn equivalent to $\llbracket \Gamma(\pi(E^\perp),\Gamma(L_\Phi)\rrbracket_{E^\perp} \subseteq \Gamma(L_\Phi)$. Therefore, by Remark \ref{rem:65} this condition is automatically satisfied since $\Phi$ is basic with respect to $E$. Thus, $J_\Phi$ is a transverse generalized complex structure if and only if $\Phi$ is a weak generalized CRF-structure.
\end{example}

\begin{theorem}\label{thm:73}
Let $E$ be a generalized almost product structure of rank $2k$ on $M$ and let $\Phi$ be a generalized $F$-structure on $E$. The following are equivalent:
\begin{enumerate}[1)]
\item $\Phi$ is a weak generalized CRF-structure on $E$;
\item $d_E(U_\Phi^k)\subseteq U_\Phi^{k-1}$;
\item $d_E\in \EE_\Phi^{-1}\oplus \EE_\Phi^{1}$;
\item $\ad_\Phi^2(d_E)=-d_E$;
\item $\llbracket \Phi,\Phi\rrbracket_E=d_E$.
\end{enumerate}
\end{theorem}

\smallskip\noindent\emph{Proof:} 1) holds if and only if $\llbracket l_1,l_2\rrbracket_E$ preserves $U_\Phi^k$ for each $l_1,l_2\in \Gamma(L_\Phi)$. Unraveling the definition of $\llbracket \,,\,\rrbracket$, 1) is equivalent to $d_E(U_\Phi^k)\subseteq \ker(\wedge^2 L_\Phi) = U^k_\Phi\oplus U_\Phi^{k-1}$. Upon inspection of $E$-biparities we conclude that 1) is equivalent to 2). It is clear that 3) is equivalent to 
\begin{equation}\label{eq:30}
d_E(U^i_\Phi)\subseteq U^{i-1}_\Phi\oplus U_\Phi^{i+1} \quad\textrm{ for all }\quad  i=0,\ldots,k 
\end{equation}
so that 2) is a particular case of 3). For the converse, assume \eqref{eq:10} holds for all $i\le j$. From the equivalence of 1) and 2) we deduce that $\llbracket l_1,l_2\rrbracket_{_E}(U^i_\Phi)\subseteq U^{i+1}_\Phi$ and thus $(l_1\circ l_2\circ d_E)(U_\Phi^i)\subseteq \Gamma(F_\Phi^{k-i-1})$ for all $l_1,l_2\in \Gamma(L_\Phi)$. This implies $d_E(U_\Phi^i)\subseteq \Gamma(F_\Phi^{k-i+1})$. Since $d_E$ is real, taking complex conjugates and inspecting $E$-biparities we conclude that \eqref{eq:10} holds and thus 2) is equivalent to 3). Assume that 3) holds and let $\partial_\Phi$ be the projection of $d_E$ onto $\EE_\Phi^1$. By Lemma \ref{lem:63}
\begin{equation}\label{eq:31}
\ad_\Phi^2(d_E)=\sqrt{-1}\ad_\Phi(\partial_\Phi-\overline\partial_\Phi)=-\partial_\Phi-\overline\partial_\Phi=-d_E\,.
\end{equation}
Conversely, assume that 4) holds and let 
\begin{equation}\label{eq:32}
\partial_\Phi = \frac{1}{2}(d_E-\sqrt{-1}\ad_\Phi(d_E))\,.
\end{equation}
From the reality of $d_E$ and $\Phi$ we obtain $d_E=\partial_\Phi+\overline\partial_\Phi$ as well as $\ad_\Phi(\partial_\Phi)=\sqrt{-1}\partial_\Phi$. A further application of Lemma \ref{lem:63} shows that 4) implies 3). The equivalence of 4) and 5) is straightforward.

\begin{example}\label{ex:74}
Let $\Phi$ be a weak generalized CRF-structure on $\T M$. Then $d=d_E=\partial_\Phi+\overline\partial_\Phi$ coincide with the decomposition of the de Rham operator induced by the generalized complex structure $J_\Phi$ given in \cite{G} and $\ad_\Phi(d)=-d^{J_\Phi}$. Moreover in this case the equivalence of 1) and 5) in Theorem \ref{thm:73} is proved in \cite{guttenberg}.   
\end{example}

\begin{example}\label{ex:75}
Let $E$ be a generalized almost product structure of rank $2$ on $M$ and let $\Phi$ be a generalized $F$-structure on $E$. Then the $\Phi$-grading is concentrated in degrees $\{0,\pm 1\}$ and since $d_E$ has by definition odd $E$-parity, it follows that from Remark \ref{rem:62} that condition 3) in Theorem \ref{thm:73} is satisfied and thus $\Phi$ is automatically a weak generalized CRF-structure. In particular, every generalized $F$-structure on a 3-manifold is an example of a weak generalized CRF-structure.
\end{example}

\begin{example}\label{ex:76}
Let $\Phi=\omega+\Lambda$ be as in Example \ref{ex:60}. It follows from $U_\Phi^k=e^{-\sqrt{-1}\omega}\Omega_{E}^{0,\bullet}$ and the equivalence of 1) and 2) in Theorem \ref{thm:73} that $\Phi$ is a weak generalized CRF-structure if and only if $d_E(e^{-\sqrt{-1}\omega}\Omega^{0,\bullet}_E)\subseteq e^{-\sqrt{-1}\omega}\Omega^{1,\bullet}_E$. Since by definition of $d_E$, $d_E(\Omega_E^{0,\bullet})\subseteq \Omega_E^{1,\bullet}$, it follows that $\Phi$ is a weak generalized CRF-structure if and only if $d_E\omega\in\Omega_E^{1,2}$. In particular, if $E$ has rank greater or equal than $2(\dim M -1)$, this condition reduces to $d_E\omega=0$. 
\end{example}

\begin{rem}\label{rem:77}
Let $\Phi$ be a weak generalized CRF-structure on an almost product structure $E$. By Lemma \ref{lem:64} and Theorem \ref{thm:73}, the de Rham operator decomposes as
\begin{equation}\label{eq:33}
d=\overline\partial_\Phi^\perp + \overline\partial_\Phi+ \delta_\Phi+\partial_\Phi+\partial_\Phi^\perp 
\end{equation}
where $\delta_\Phi\in \EE_\Phi^0$ is real, $\partial_\Phi\in \EE_\Phi^1$ and $\partial_\Phi^\perp\in \EE_\Phi^2$. Analyzing the $\Phi$ graded components of the identity $[d,d]=0$ we obtain
\begin{align}
[\partial_\Phi^\perp,\partial_\Phi^\perp]=[\partial_\Phi^\perp,\partial_\Phi] & = 0\,;\label{eq:34}\\
2[\partial_\Phi^\perp, \delta_\Phi]+[\partial_\Phi,\partial_\Phi]&=0\,;\label{eq:35}\\
[\partial_\Phi^\perp,\overline \partial_\Phi]+[\partial_\Phi,\delta_\Phi] & =0\,;\label{eq:36}\\
2[\partial_\Phi^\perp,\overline\partial_\Phi^\perp]+2[\partial_\Phi,\overline \partial_\Phi]+[\delta_\Phi,\delta_\Phi] &=0\,.\label{eq:37}
\end{align}
In particular, $(\Omega_M\otimes \C,\partial_\Phi^\perp)$ is a chain complex on which $\partial_\Phi$ acts by chain maps.
\end{rem}

\begin{example}\label{ex:78}
Let $\Phi=\omega+\Lambda$ be as in Example \ref{ex:60} and assume that $d_E\omega=0$. Condition 5) in Theorem \ref{thm:73} ensures that $\ad_\Lambda^2(d_E)=0$ and thus
\begin{equation}\label{eq:38}
\left(\Ad_{\frac{\sqrt{-1}}{2}\Lambda}\circ\Ad_{-\sqrt{-1}\omega}\right)(d_E)=\Ad_{\frac{\sqrt{-1}}{2}\Lambda}(d_E)=d_E+\frac{\sqrt{-1}}{2}\ad_\Lambda(d_E)
\end{equation}
or, equivalently,
\begin{equation}\label{eq:39}
[d_E,\Psi]=\frac{\sqrt{-1}}{2}(\Psi\circ\ad_\Lambda)(d_E)\,.
\end{equation}
Comparing $\Phi$-gradings we further obtain
\begin{equation}\label{eq:40}
\overline \partial_\Phi\circ \Psi = \Psi \circ d_E \quad \text{ and } \quad \partial_\Phi\circ \Psi =\frac{\sqrt{-1}}{2}\Psi \circ \ad_\Lambda(d_E)\,.
\end{equation}
In the case where $\omega$ is symplectic, we recover Theorem 2.3 in \cite{cavalcanti}.
\end{example}

\begin{mydef}\label{def:79}
A generalized $F$-structure $\Phi$ on a split structure $E$ on $M$ is a {\it generalized CRF-structure} if $\llbracket l_1,l_2\rrbracket\in \Gamma(L_\Phi)$ for all $l_1,l_2\in \Gamma(L_\Phi)$.
\end{mydef}

\begin{example}\label{ex:80}
Let $E$ be one of the generalized almost product structures on $M=S^3$ defined in Example \ref{ex:34} and let $\Phi$ be a generalized $F$-structure on $E$. A direct calculation shows that there exists $\tau\in \Omega_M\otimes\C\setminus\Omega_M\otimes\R$ such that either
\begin{equation}\label{eq:41}
L_\Phi={\rm span}\{\alpha_2+\tau\alpha_3,x_3-\tau x_2\}\,
\end{equation}
in which case $U^1_\Phi={\rm span}\{\alpha_2+\tau\alpha_3,\alpha_1\alpha_2+\tau\alpha_1\alpha_3\}$ has odd $E$-parity, or
\begin{equation}\label{eq:42}
L_\Phi={\rm span}\{x_2+\tau \alpha_3,x_3-\tau \alpha_2\}
\end{equation}
in which case $U^1_\Phi={\rm span}\{1+\tau\alpha_2\alpha_3,\alpha_1+\tau\alpha_1\alpha_2\alpha_3\}$ has even $E$-parity. As pointed out in Example \ref{ex:75}, in either case $L_\Phi$ is automatically involutive with respect to $\llbracket\,,\,\rrbracket_E$. Therefore, if \eqref{eq:41} holds then $\Phi$ is a generalized CRF-structure if and only if 
\begin{equation}\label{eq:43}
0=\llbracket x_3-\tau x_2,\alpha_2+\tau\alpha_3\rrbracket_{E^\perp}=2(1+\tau^2-\L_{x_1}(\tau))\alpha_1\,.
\end{equation}
In particular if $\tau$ is constant, then it must equal to $\pm \sqrt{-1}$.
On the other hand, since $[\llbracket x_2,x_3\rrbracket,\alpha_1]=2$, then
$\llbracket x_2+\tau \alpha_3,x_3-\tau \alpha_2\rrbracket_{E^\perp}\neq 0$ and thus $\Phi$ is never a generalized CRF-structure if \eqref{eq:42} holds.
\end{example}

\begin{theorem}\label{thm:81}
Let $E$ be a generalized almost product structure of rank $2k$ on $M$ and let $\Phi$ be a generalized $F$-structure on $E$. The following are equivalent:
\begin{enumerate}[1)]
\item $\Phi$ is a generalized CRF-structure on $E$;
\item $d(U_\Phi^k)\subseteq U_\Phi^k\oplus U_\Phi^{k-1}$;
\item $d\in \EE_\Phi^{-1}\oplus \EE_\Phi^0\oplus \EE_\Phi^{1}$;
\item $(\ad_\Phi^3+\ad_\Phi)(d)=0$;
\item $\llbracket \Phi,\Phi\rrbracket = d_E$.
\end{enumerate}
\end{theorem}

\smallskip\noindent\emph{Proof:} Since $\Phi$ is a generalized CRF-structure on $E$ if and only if $\Phi$ is a weak generalized CRF-structure on $E$ such that $\llbracket l_1,l_2\rrbracket_{E^\perp}=0$ for all $l_1,l_2\in \Gamma(L_\Phi)$, it follows from Theorem \ref{thm:73} and Lemma \ref{lem:64} that the first three statements are equivalent and that any of them implies 4). Conversely, assume that 4) holds and set
\begin{equation}\label{eq:44}
 \partial_\Phi=\frac{1}{2}({\rm ad}^2_\Phi+\sqrt{-1}{\rm ad}_\Phi)(d)\,.
\end{equation}
By Lemma \ref{lem:63}, $\partial_\Phi\in \EE_\Phi^1$ and thus $\overline\partial_\Phi\in \EE_\Phi^{-1}$. A further application of Lemma \ref{lem:63} yields $d-\partial_\Phi-\overline\partial_\Phi=d+\ad_\Phi^2(d)\in \EE_\Phi^0$, from which we conclude that 4) implies 3). Inspecting $E$-biparities and using Theorem \ref{thm:73} shows that 5) is equivalent to the statement that $\Phi$ is a weak generalized CRF-structure and, using Lemma \ref{lem:64},
\begin{equation}\label{eq:45}
8{\rm Re}(\partial^\perp_\Phi)=-\ad_\Phi^2(d_{E^\perp})=\llbracket \Phi,\Phi \rrbracket_{E^\perp}=0\,.
\end{equation}
Since $d_{E^\perp}-\delta_\Phi=2{\rm Re}(\partial_\Phi^\perp)$, \eqref{eq:45} shows that 3) is equivalent to 5) and the Theorem is proved.

\begin{rem}\label{rem:82}
Let $E$ be a generalized almost product structure on $M$ and let $\Phi$ be a generalized CRF-structure on $E$. Then $d_{E^\perp}=\delta_\Phi$ and  $d_E=\partial_\Phi +\overline \partial_\Phi = \ad_\Phi^2(d)$.  From \eqref{eq:35} and \eqref{eq:36} we conclude that $(\Omega_M,\partial_\Phi)$ is a complex on which $d_{E^\perp}$ acts by chain maps. 
\end{rem}

\begin{example}\label{eq:83}
Let $\Phi=\omega+\Lambda$ be as in Example \ref{ex:60}. By Theorem $\Phi$ is a generalized CRF-structure if and only if $d(e^{-\sqrt{-1}\omega}\Omega^{0,\bullet}_E)\subseteq e^{-\sqrt{-1}\omega}\Omega^{1,\bullet}_E$. Using Example \ref{ex:76}, this is equivalent to the conditions 
\begin{enumerate}[i)]
\item $d_E\omega\in \Omega^{1,2}_E$;
\item $d_{E^\perp}\omega=0$;
\item $d_{E^\perp}\Omega_E^{0,\bullet}\subseteq \Omega_E^{0,\bullet}$ i.e. $\pi(E)$ is a foliation.  
\end{enumerate}
If $\Lambda\in \Gamma(\wedge^2 TM)$, then condition 5) in Theorem \ref{thm:81} implies that $\Lambda$ is Poisson. 
\end{example}

\begin{rem}\label{rem:84}
Let $\Phi=\omega+\Lambda$ be as in Example \ref{ex:60} be such that $\pi(E)$ is a foliation and $d\omega=0$. Arguing as in Example \ref{ex:78}, $[d,\Psi]=\frac{\sqrt{-1}}{2}\ad_\Lambda(d)$. Matching $\Phi$-grading we obtain, in addition to \eqref{eq:39}, $\ad_\Lambda(d_E^\perp)=0$ and thus $[d_{E^\perp},\Psi]=0$. If $\Lambda$ is a bivector, then $(\Omega_M,\ad_{\Lambda(d)})$ is the complex that computes the canonical homology of the Poisson manifold $(M,\Lambda)$ defined in \cite{brylinski}. Since complex conjugation is an isomorphism, we conclude that the canonical homology of $(M,\Lambda)$ is isomorphic to the cohomology of the complex $(\Omega_M,d_E)$. In the symplectic case, we obtain the isomorphism with de Rham cohomology noticed in \cite{brylinski}. In the cosymplectic case, this is proved in \cite{fernandez98}.  
\end{rem}

\begin{rem}\label{rem:85}
Let $\Phi$ be a generalized $F$-structure on a split structure $E$ and fix $D_1,D_2\subseteq E^\perp$ as in Remark \ref{rem:67}. It is easy to adapt the arguments of this section to prove that the equivalence of conditions 1)-4) in Theorem \ref{thm:81} holds in this more general setting. Notice that the notion of generalized CRF-structure is invariant under T-duality.
\end{rem}

\begin{prop}\label{prop:86}
Let be $E$ be a generalized almost product structure on $M$ and let $\Phi$ be a generalized CRF-structure on $E$. The following are equivalent:
\begin{enumerate}[1)]
\item $E$ is a generalized local product structure;
\item $[\partial_\Phi,\overline\partial_\Phi]=0$;
\item $[\ad_\Phi(d),\ad_\Phi(d)]=0$.
\end{enumerate} 
\end{prop}

\smallskip\noindent\emph{Proof:} The equivalence of 1) and 2) follows from Proposition \ref{prop:45} and \eqref{eq:37}. Since $\Phi$ is generalized CRF, using \eqref{eq:44} we obtain
\begin{equation}\label{eq:47}
4[\partial_\Phi,\overline\partial_\Phi]=[\ad_\Phi^2(d),\ad_\Phi^2(d)]+[\ad_\Phi(d),\ad_\Phi(d)]=[d_E,d_E]+[\ad_\Phi(d),\ad_\Phi(d)]\,.
\end{equation}
Plugging in \eqref{eq:37}, we obtain $6[\partial_\Phi,\overline\partial_\Phi]=[\ad_\Phi(d),\ad_\Phi(d)]$ which concludes the proof.

\begin{rem}\label{rem:87}
Let $\Phi$ be a generalized CRF-structure on a generalized local product structure $E$ on $M$. Consider the {\it periodic bicomplex} $(P_\Phi^{\bullet,\bullet},\partial_\Phi,\overline\partial_\Phi)$ whose bigraded components are $P_\Phi^{p,q}=U_\Phi^{p-q}$. Since $P_\Phi^{\bullet,\bullet}$ is bounded in both directions, the corresponding spectral sequence converges to the cohomology of $d_E$. In the case $E=\T M$ we recover the canonical spectral sequence of \cite{cavalcanti}. In fact it is easy to show that the considerations of Sections 4 and 5 in \cite{cavalcanti} extend verbatim to this more general setting and so we conclude that the spectral sequence of the periodic bicomplex degenerates at the first page if the {\it $\partial_\Phi\overline\partial_\Phi$-lemma} holds i.e.\ if
\begin{equation}\label{eq:48}
{\rm Im}(\partial_\Phi)\cap\ker(\overline \partial_\Phi)={\rm Im}(\overline\partial_\Phi)\cap\ker(\partial_\Phi)={\rm Im}(\partial_\Phi\overline\partial_\Phi)
\end{equation}
or equivalently if the inclusion of complexes $(\Omega_M^\bullet\cap \ker(\ad_\Phi(d_E)),d_E)\hookrightarrow(\Omega_M^\bullet, d_E)$ is a quasi-isomorphism. Conversely if the spectral sequence of the periodic bicomplex degenerates at the first page and the $\Phi$-grading induces a splitting of cohomology, then the $\partial_\Phi\overline\partial_\Phi$-lemma holds.
\end{rem}

\begin{example}\label{ex:88}
Let $\Phi=\omega+\Lambda$ be a generalized CRF-structure as in Remark \ref{rem:84}. Since $d\omega=0$, then then $E$ is a generalized local product structure. Moreover, by \eqref{eq:40}, $(U_\Phi^\bullet,\overline\partial_\Phi)$ is isomorphic to $(\Omega_M, d_E)$. Therefore, the spectral sequence of the periodic bicomplex degenerates at the first page (even though the $\partial_\Phi\overline\partial_\Phi$-lemma does not hold in general).  
\end{example}

\begin{rem}\label{rem:89}
Let $\Phi$ be a weak generalized CRF-structure on a generalized almost product structure $E$ such that $\pi(E^\perp)$ is a foliation. If $\Phi$ is basic with respect to $E$, one can repeat the construction of Remark \ref{rem:87} and define the {\it basic periodic bicomplex} $(\mathcal BP^{\bullet,\bullet}_\Phi,\partial_\Phi,\overline\partial_\Phi)$, where $\mathcal BP^{p,q}_\Phi=\mathcal B_\Phi^{p-q}$. Then once again the calculations of \cite{cavalcanti} apply and one concludes that \eqref{eq:48} holds for the restrictions of $\partial_\Phi$ and $\overline\partial_\Phi$ to $\mathcal B_E$ if and only if the spectral sequence of the basic periodic bicomplex degenerates at the first page and the $\Phi$-grading induces a cohomological grading on $H(\mathcal B_E,d_E)$. In the case $B_E=0$, this is the main result of \cite{razny16}.
\end{rem}

\begin{rem}\label{rem:90}
Let $p:M\to N$ be a fiber bundle with Ehresmann connection $H\subseteq TM$ and let $E$ be the almost product structure on $M$ generated by $p^*\Omega^1_N$ and $\Gamma(H)$. Any generalized almost complex structure $\Psi$ on $N$ defines a generalized $F$-structure $\Phi$ on $E$ by setting $\Phi\circ p^*=p^*\circ \Psi$ and $\Phi({\rm Ann}(H))=0$. Using $d_E\circ p^*=d\circ p^*=p^*\circ d$, we obtain
\begin{equation}\label{eq:49}
\llbracket\Phi,\Phi\rrbracket \circ p^*=\llbracket \Phi,\Phi\rrbracket_E \circ p^*=p^*\circ \llbracket \Psi,\Psi\rrbracket\,.  \end{equation}
Since $p^*:\Omega_N\to\mathcal B_E$ is an isomorphism, we conclude that if $\Phi$ is a weak generalized CRF-structure on $E$, then $\Psi$ is a generalized complex structure on $N$. Conversely, if $\Psi$ is a generalized complex structure on $N$, then \eqref{eq:49} shows that $(\llbracket \Phi,\Phi\rrbracket_E-d_E)(\mathcal B_E)=0$. On the other  hand, $d_E(\Omega^{0,1})\subseteq \Omega^{1,1}$, then $\ad_\Phi^2(d_E)+d_E$ vanishes on $\Omega^{0,1}$. Since $\Omega_M$ is generated by $\mathcal B_E\otimes \Omega^{0,\bullet}$, we conclude by Theorem \ref{thm:73} that $\Psi$ is a generalized complex structure on $N$ if and only if $\Phi$ is a weak generalized CRF-structure on $E$.   
\end{rem}

\begin{prop}\label{prop:91}
Let $E$ be a generalized local product structure on $M$ such that $dB_E=0$ and let $\Phi$ be a generalized $F$-structure on $E$. The following are equivalent:
\begin{enumerate}[1)]
\item $\Phi$ is a generalized CRF-structure;
\item $\Phi$ is a weak generalized CRF-structure;
\item $\Ad_{B_E}(\Phi)$ induces a generalized complex structure on the leaves of $\pi(E)$ and acts trivially on the leaves of $\pi(E^\perp)$. 
\end{enumerate}
\end{prop}

\smallskip\noindent\emph{Proof:} Clearly, 1) implies 2). Furthermore, using \eqref{eq:9}, Theorem \ref{thm:73} and the assumption $dB_E=0$ we conclude that $\Phi$ is a weak generalized CRF-structure on $E$ if and only if $\Ad_{B_E}(\Phi)$ is a weak generalized CRF-structure on $\Ad_{B_E}(E)$. Therefore we may assume $B_E=0$ so that the equivalence of 2) and 3) is given by Remark \ref{rem:90}. On the other hand, if condition 3) holds, then \eqref{eq:49} implies $(\llbracket \Phi,\Phi\rrbracket -d_E)(\mathcal B_E)=0$. Moreover, $(\llbracket \Phi,\Phi\rrbracket -d_E)$ is $\Omega_M^0$-linear and vanishes on $\mathcal B_{E^\perp}$. Since $\Omega_M$ is locally generated by $\mathcal B_E\otimes \mathcal B_{E^\perp}$ we conclude by Theorem \ref{thm:81} that $\Phi$ is a generalized CRF-structure.

\begin{theorem}\label{thm:92}
Let $\Phi$ be a generalized $F$-structure on $\T M$ such that $J_\Phi$ is a generalized complex structure. Then $(M,J_\Phi)$ is (possibly up to a $B$-field transform by a closed $2$-form) locally the product of two generalized complex manifolds if and only if there exists a generalized almost product structure $E$ on $M$ and a generalized CRF-structure $\Phi_E$ on $E$ such that 
\begin{enumerate}[1)]
\item $dB_E=0$;
\item $(J_\Phi-J_{\Phi_E})(E)=0$;
\item $[\partial_{\Phi_E},\overline \partial_{\Phi_E}]$.
\end{enumerate}
\end{theorem}

\smallskip\noindent\emph{Proof:} Suppose that $J_\Phi$ restricts to generalized complex structures on the leaves of two complementary foliations. We may assume that the foliations are of the form $\pi(E),\pi(E^\perp)$, for some generalized local product structure $E$. the corresponding. By construction, 1) holds. By Lemma \ref{lem:53} there exist a generalized $F$-structure $\Phi_E$ on $E$ such that 2) holds. By Proposion \ref{prop:91}, $\Phi_E$ is a generalized CRF-structure. Finally, 3) holds by Proposition \ref{prop:86}. Using \eqref{eq:9}, the same conditions hold if a $B$-field transform by a closed $2$-form is applied to $J_\Phi$. Conversely, suppose that $E$ is a generalized almost product structure on $M$ and $\Phi_E$ is a generalized CRF-structure on $E$ such that conditions 1)-3) hold. By Proposition \ref{prop:86} it follows that $E$ is a local product structure and by Proposition \ref{prop:91} $\Ad_{B_E}(\Phi_E)$ induces generalized complex structures on the leaves of $\pi(E)$ (while acting trivially on the leaves of $\pi(E^\perp)$). By Theorem \ref{thm:81}, $[\Phi_E,d_{E^\perp}]=0$ and using $[\Phi_E,\Phi-\Phi_E]=0$ (together with \eqref{eq:5}) we obtain 
\begin{equation}\label{eq:50}
d_{E^\perp}=\llbracket \Phi,\Phi\rrbracket_{E^\perp}=\llbracket \Phi-\Phi_E,\Phi-\Phi_E\rrbracket_{E^\perp}\,. 
\end{equation}
Therefore, $\Phi-\Phi_E$ is a weak generalized CRF-structure on $E^\perp$. By Proposition \ref{prop:91}, $\Ad_{B_E}(\Phi_E)$ induces a generalized complex structure on the leaves of $\pi(E)$ and acts trivially on the leaves of $\pi(E^\perp)$. Thus the generalized complex structure $J_{\Ad_{B_E}(\Phi)})$ is locally the product of the generalized complex structures $J_{\Ad_{B_E}(\Phi_E)}$ and $J_{\Ad_{B_E}(\Phi-\Phi_E)}$, on the leaves of $\pi(E)$ and $\pi(E^\perp)$, respectively.

\section{Deformations}

\begin{rem}\label{rem:93}
Let $\Phi_1,\Phi_2$ be a generalized $F$-structures on an almost product structure $E$ on $M$ and  For $\Phi_2-\Phi_1$ ``small'' there exists $\epsilon \in \Gamma(\wedge^2 \overline L_\Phi)$ such that $\Gamma(L_{\Phi_2})=\Ad_\epsilon(\Gamma(L_{\Phi_1}))$. If this is the case, then $U_{\Phi_2}^i=e^\epsilon(U_{\Phi_1}^i)$ for all $i$.
\end{rem}

\begin{lem}\label{lem:94}
Let $\Phi$ be a generalized $F$-structure on an almost product structure $E$ on $M$ and let $\alpha,\beta\in \Gamma(\wedge^\bullet\overline L_\Phi)$. Then 
\begin{enumerate}[1)]
\item $\llbracket \alpha,\beta\rrbracket_{E^\perp}=\llbracket \alpha,\beta\rrbracket_{\partial_\Phi^\perp}$;
\item $\Ad_\alpha(d_{E^\perp})=d_{E^\perp}+\ad_\alpha(d_{E^\perp})+\frac{1}{2}\ad_\alpha^2(\partial_\Phi^\perp)$.
\end{enumerate}
\end{lem}

\smallskip\noindent\emph{Proof:} Lemma \ref{lem:64} implies the decomposition $d_{E^\perp}=\overline \partial_\Phi^\perp+\delta_\Phi+\partial_\Phi^\perp$ in to $\Phi$-graded components. Since $\L^{\!\!^{E^\perp}}_\varphi$ is a derivation for every operator $\varphi$ and \eqref{eq:5} holds, it suffices to establish 1) if $\alpha,\beta\in \Gamma(\overline L_\Phi)$. In this case, $\llbracket \alpha,\beta\rrbracket_{E^\perp}\in \Gamma(E^\perp)\otimes \C\subseteq \EE_\Phi^0$ which proves 1). The second statement follows from the first since
\begin{equation}\label{eq:51}
\ad_\alpha^2(d_{E^\perp})=-\llbracket\alpha,\alpha\rrbracket_{E^\perp}=\ad_\alpha^2(\partial_\Phi^\perp)
\end{equation}
is a section of $E^\perp\otimes\C$ and thus commutes with $\alpha$.

\begin{lem}\label{lem:95}
Let $E$ be a generalized almost product structure on $M$, let $\Phi$ be a generalized CRF-structure and let $\alpha,\beta\in \Gamma(\wedge^\bullet \overline L_\Phi)$. Then
\begin{enumerate}[1)]
\item $\llbracket \alpha,\beta\rrbracket_E = \llbracket \alpha,\beta\rrbracket_{\partial_\Phi}$;
\item $\Ad_\alpha(d_E)=d_E+\ad_\alpha(d_E)+\frac{1}{2}\ad_\alpha^2(\partial_\Phi)$.
\end{enumerate}
\end{lem}

\smallskip\noindent\emph{Proof:} By linearity it suffices to consider the case $\alpha\in \Gamma(\wedge^a \overline L_\Phi)$ and $\beta\in \Gamma(\wedge^b\overline L_\Phi)$, with $a$ and $b$ arbitrary non-negative integers. Since $\Phi$ is a weak generalized CRF-structure, then $\llbracket \alpha,\beta\rrbracket_E\in \Gamma(\wedge^{a+b-1}\overline L_\Phi)$. Taking into account the decomposition \eqref{eq:33} this proves 1). The second statement follows from the first since $\ad_\alpha^2(d_E)=-\llbracket \alpha,\alpha\rrbracket_E=\ad_\alpha^2(\partial_\Phi)$ is a section of $\wedge^\bullet \overline L_\Phi$.

\begin{prop}\label{prop:96}
Let $E$ be a generalized almost product structure on $M$, let $\Phi$ be a weak generalized CRF-structure and let $\Phi+\Phi_\epsilon$ be a generalized $F$-structure on $E$ such that $\Gamma(L_{\Phi+\Phi_\epsilon})=\Ad_\epsilon(\Gamma(L_\Phi))$ for some $\epsilon\in \Gamma(\wedge^2\overline L_\Phi)$. The following are equivalent:
\begin{enumerate}[1)]
\item $\Phi+\Phi_\epsilon$ is a weak generalized CRF-structure;
\item $[\epsilon,\overline\partial_\Phi]+\frac{1}{2}\llbracket \epsilon,\epsilon\rrbracket_E=0$;
\item $[\Phi_\epsilon,[\Phi,d_E]]+\frac{1}{2}[[\Phi,\Phi_\epsilon],d_E]+\frac{1}{2}[\Phi_\epsilon,[\Phi_\epsilon,d_E]]=0$;
\item $\llbracket\Phi,\Phi_\epsilon \rrbracket_E + \llbracket \Phi_\epsilon,\Phi\rrbracket_E+\llbracket \Phi_\epsilon,\Phi_\epsilon\rrbracket_E =0$.
\end{enumerate}
\end{prop}

\smallskip\noindent\emph{Proof:} By Theorem \ref{thm:73}, $\Phi+\Phi_\epsilon$ is a weak generalized CRF-structure if and only if $d_E (e^\epsilon U^i_\Phi)\subseteq e^\epsilon (U_\Phi^{i+1}\oplus U^{i-1}_\Phi)$ for all $i$ if and only if $\Ad_{-\epsilon}(d_E)\subseteq \EE_\Phi^{-1}\oplus \EE_\Phi^1$. By Lemma \ref{lem:95}, this is equivalent to the vanishing of the term of $\Phi$-degree $-3$ in $\Ad_{-\epsilon}(d_E)$ i.e.\ the expression to the LHS of 2). Therefore 1) and 2) are equivalent. Using again Theorem \ref{thm:73}, 1) is also equivalent to
\begin{equation}\label{eq:52}
\llbracket \Phi,\Phi\rrbracket_E+\llbracket\Phi,\Phi_\epsilon \rrbracket_E + \llbracket \Phi_\epsilon,\Phi\rrbracket_E+\llbracket \Phi_\epsilon,\Phi_\epsilon\rrbracket_E=\llbracket\Phi+\Phi_\epsilon,\Phi+\Phi_\epsilon\rrbracket_E=d_E\,.
\end{equation}
Since the assumption that $\Phi$ is a weak generalized CRF-structure implies $\llbracket \Phi,\Phi\rrbracket_E=d_E$, \eqref{eq:52} is equivalent to 4). A straightforward calculation involving \eqref{eq:3}, shows that 3) is equivalent to 4) and the Theorem is proved.

\begin{example}\label{ex:97}
If $\Phi$ be weak generalized CRF-structure on $\T M$. The equivalent conditions of Proposition \ref{prop:96} coincide with the various forms of the Kodaira-Spencer equation for the generalized complex structure $J_\Phi$ given in \cite{tomasiello}. 
\end{example}

\begin{rem}\label{rem:98}
Let $E$ be a generalized almost product structure on $M$ such that $\pi(E)$ is a foliation and let $\Phi$ be a weak generalized CRF-structure on $E$. By Proposition \ref{prop:41}, $E$ is a Courant algebroid with respect to $\llbracket\,,\,\rrbracket_E$. As shown in \cite{LWX}, it follows that $(L_\Phi,\overline L_\Phi)$ is a Lie bialgebroid and infinitesimal deformations are solutions to the Maurer-Cartan equation 
\begin{equation}\label{eq:53}
d_{L_\Phi}(\epsilon)+\frac{1}{2}\llbracket \epsilon,\epsilon \rrbracket_E=0
\end{equation}
where $d_{L_\Phi}$ is the Lie algebroid differential on $\wedge^\bullet \overline L_\Phi$ obtained by identifying $\overline L_\Phi$ with $L_\Phi^*$. This is compatible with Proposition \ref{prop:96} since
\begin{align}
[l_3,[l_2,[l_1,[\epsilon,\overline \partial_\Phi]]]]&=[l_3,[\llbracket l_1,l_2\rrbracket_E,\epsilon]]-[l_2,[\llbracket l_1,l_3\rrbracket_E,\epsilon]]+[l_1,[\llbracket l_2,l_3\rrbracket_E,\epsilon]]\nonumber\\
& - \llbracket l_1,[l_3,[l_2,\epsilon]]\rrbracket_E +\llbracket l_2,[l_3,[l_1,\epsilon]]\rrbracket_E - \llbracket l_3,[l_2,[l_1,\epsilon]]\rrbracket_E\label{eq:54}
\end{align}
for all $l_1,l_2,l_3\in \Gamma(L_\Phi)$.
\end{rem}

\begin{rem}\label{rem:99}
Let $E$ be a generalized almost product structure such that $\pi(E^\perp)$ is a foliation and let $\Phi$ be a weak generalized CRF-structure on $E$ that is basic with respect to $E$. By Proposition \ref{prop:96} infinitesimal deformations of transverse generalized complex structures are parametrized by operators $\epsilon\in \Gamma(\wedge^2\overline L_\Phi)$ that are basic with respect to $E$ and satisfy $[\epsilon,\overline \partial_\Phi]+\frac{1}{2}\llbracket \epsilon,\epsilon \rrbracket_E=0$.
\end{rem}

\begin{theorem}\label{thm:100}
Let $E$ be a generalized almost product structure on $M$, let $\Phi$ be a generalized CRF-structure on $E$ and let $\Phi+\Phi_\epsilon$ be a generalized $F$-structure on $E$ such that $\Gamma(L_{\Phi+\Phi_\epsilon})=\Ad_\epsilon(\Gamma(L_\Phi))$ for some $\epsilon\in \Gamma(\wedge^2\overline L_\Phi)$. The following are equivalent:
\begin{enumerate}[1)]
\item $\Phi+\Phi_\epsilon$ is a generalized CRF-structure;
\item $[\epsilon,\overline \partial_\Phi]+\frac{1}{2}\llbracket \epsilon,\epsilon \rrbracket=0$ and $[\epsilon,d_{E^\perp}]=0$;
\item $[\Phi_\epsilon,[\Phi,d]]+\frac{1}{2}[[\Phi,\Phi_\epsilon],d]+\frac{1}{2}[\Phi_\epsilon,[\Phi_\epsilon,d]]=0$;
\item $\llbracket\Phi,\Phi_\epsilon \rrbracket + \llbracket \Phi_\epsilon,\Phi\rrbracket+\llbracket \Phi_\epsilon,\Phi_\epsilon\rrbracket =0$.
\end{enumerate}
\end{theorem}

\smallskip\noindent\emph{Proof:} Since $\Gamma(\overline L_\Phi)$ is closed under the Dorfman bracket, then $\llbracket \epsilon,\epsilon\rrbracket_E=\llbracket \epsilon,\epsilon\rrbracket$ and thus $\Phi+\Phi_\epsilon$ is a weak generalized CRF-structure if and only if $[\epsilon,\overline\partial_\Phi]+\frac{1}{2}\llbracket \epsilon,\epsilon\rrbracket=0$ by Proposition \ref{prop:96}. Moreover, since $\Phi$ is a generalized CRF-structure, then $d_{E^\perp}\in \EE_{\Phi+\Phi_\epsilon}^0$ if and only if $\Ad_{-\epsilon}(d_{E^\perp})\in \EE_\Phi^0$. By Lemma \ref{lem:94}, this occurs if and only if $[\epsilon,d_{E^\perp}]=0$, which proves the equivalence of 1) and 2). The equivalence of 4) and 1) is straightforward from Theorem \ref{thm:81} while the equivalence of 4) and 5) follows from a direct calculation involving \eqref{eq:3}.


\begin{bibdiv} 
\begin{biblist}

\bib{AG}{article}{
   author={Aldi, Marco},
   author={Grandini, Daniele},
   title={Generalized contact geometry and T-duality},
   journal={J. Geom. Phys.},
   volume={92},
   date={2015},
   pages={78--93}
}

\bib{AG2}{article}{
   author={Aldi, Marco},
   author={Grandini, Daniele},
   title={An abstract Morimoto theorem for generalized $F$-structures},
   journal={Q. J. Math.},
   volume={67},
   date={2016},
   number={2},
   pages={161--182}
}




\bib{brylinski}{article}{
   author={Brylinski, Jean-Luc},
   title={A differential complex for Poisson manifolds},
   journal={J. Differential Geom.},
   volume={28},
   date={1988},
   number={1},
   pages={93--114},
   issn={0022-040X},
   review={\MR{950556}},
}


\bib{cavalcanti}{article}{
   author={Cavalcanti, Gil R.},
   title={The decomposition of forms and cohomology of generalized complex
   manifolds},
   journal={J. Geom. Phys.},
   volume={57},
   date={2006},
   number={1},
   pages={121--132}
}



\bib{fernandez98}{article}{
   author={Fern{\'a}ndez, Marisa},
   author={Ib{\'a}{\~n}ez, Ra{\'u}l},
   author={de Le{\'o}n, Manuel},
   title={The canonical spectral sequences for Poisson manifolds},
   journal={Israel J. Math.},
   volume={106},
   date={1998},
   pages={133--155}
}



\bib{G}{article}{
author={Gualtieri, Marco},
title={Generalized complex geometry},
journal={Ann. of Math. (2)},
volume={174},
date={2011},
number={1},
pages={75--123}
}

\bib{gugenheim-spencer}{article}{
   author={Gugenheim, V. K. A. M.},
   author={Spencer, D. C.},
   title={Chain homotopy and the de Rham theory},
   journal={Proc. Amer. Math. Soc.},
   volume={7},
   date={1956},
   pages={144--152}
}

\bib{guttenberg}{article}{
   author={Guttenberg, Sebastian},
   title={Brackets, sigma models and integrability of generalized complex
   structures},
   journal={J. High Energy Phys.},
   date={2007},
   number={6},
   pages={004, 67 pp. (electronic)}
}

\bib{hitchin03}{article}{
   author={Hitchin, Nigel},
   title={Generalized Calabi-Yau manifolds},
   journal={Q. J. Math.},
   volume={54},
   date={2003},
   number={3},
   pages={281--308}
}

\bib{li05}{article}{
	author={Yi,Li},
    title={Deformations of Generalized Complex Structures: the Calabi-Yau case},
    eprint={arXiv:hep-th/0508030}
}



\bib{LWX}{article}{
   author={Liu, Zhang-Ju},
   author={Weinstein, Alan},
   author={Xu, Ping},
   title={Manin triples for Lie bialgebroids},
   journal={J. Differential Geom.},
   volume={45},
   date={1997},
   number={3},
   pages={547--574}
}



\bib{PW}{article}{
   author={Poon, Yat Sun},
   author={Wade, A{\"{\i}}ssa},
   title={Generalized contact structures},
   journal={J. Lond. Math. Soc. (2)},
   volume={83},
   date={2011},
   number={2},
   pages={333--352}
}

\bib{razny16}{article}{
	author={Ra\'zny, Pawel},
    title={The basic $dd^{\mathcal J}$-lemma},
    eprint={arXiv:1609.04539}
}


\bib{reinhart}{article}{
   author={Reinhart, Bruce L.},
   title={Harmonic integrals on foliated manifolds},
   journal={Amer. J. Math.},
   volume={81},
   date={1959},
   pages={529--536}
}


\bib{tomasiello}{article}{
   author={Tomasiello, Alessandro},
   title={Reformulating supersymmetry with a generalized Dolbeault operator},
   journal={J. High Energy Phys.},
   date={2008},
   number={2},
   pages={010, 25}
}





\bib{vaismanCRF}{article}{
   author={Vaisman, Izu},
   title={Generalized CRF-structures},
   journal={Geom. Dedicata},
   volume={133},
   date={2008},
   pages={129--154}
}



\bib{wade}{article}{
   author={Wade, A{\"{\i}}ssa},
   title={Local structure of generalized contact manifolds},
   journal={Differential Geom. Appl.},
   volume={30},
   date={2012},
   number={1},
   pages={124--135}
}

\bib{wade10}{article}{
   author={Wade, A{\"{\i}}ssa},
   title={Transverse geometry and generalized complex structures},
   journal={Afr. Diaspora J. Math. (N.S.)},
   volume={9},
   date={2010},
   number={2},
   pages={139--149}
}

\end{biblist}
\end{bibdiv}
\end{document}